\documentclass[a4paper,leqno,10pt]{article}
\usepackage{amsmath,amsfonts,amssymb}
\usepackage{fancyhdr}
\usepackage{makeidx}
\makeindex 

\usepackage{xypic}
\newcommand{\X}{{\mathcal X}}
\newcommand{\Hi}{{\mathcal H}}

\newcommand{\Cl}{{\mathcal C}}
\newcommand{\RR}{{\mathcal R}}
\newcommand{\LL}{{\mathcal L}}

\DeclareMathOperator{\OO}{\mathcal O}
\DeclareMathOperator{\FF}{\mathcal F}
\DeclareMathOperator{\GG}{\mathcal G}
\DeclareMathOperator{\KK}{\mathcal K}
\DeclareMathOperator{\Ss}{\mathcal S}

\DeclareMathOperator{\Pp}{\mathcal P}
\DeclareMathOperator{\DD}{\mathcal D}

\newcommand{\ta}{{}^t\!}

\newcommand{\re}{representation }
\newcommand{\res}{representations }

\newcommand{\pr}{{\bf Proof}}

\newcommand{\GL}{{\rm GL}}
\newcommand{\M}{{\rm M}}

\newcommand{\Heis}{{\rm Heis} }
\newcommand{\SL}{{\rm SL}}
\newcommand{\Spe}{{\rm Sp}}

\newcommand{\sgn}{{\rm sgn}}

\newcommand{\Ad}{{\rm Ad}}

\newcommand{\Lie}{{\rm Lie\,}}
\newcommand{\Sym}{{\rm Sym}}

\newcommand{\Tr}{{\rm Tr } }

\newcommand{\Q}{{\bf Q}}

\newcommand{\R}{{\bf R}}
\newcommand{\C}{{\bf C}}
\newcommand{\Z}{{\bf Z}}
\newcommand{\N}{{\bf N}}

\newcommand{\Ha}{{\mathfrak H}}
\newcommand{\D}{{\mathfrak D}}

\newcommand{\Or}{{\rm O}}
\newcommand{\Un}{{\rm U}}
\newcommand{\Mp}{{\rm Mp}}
\newcommand{\SO}{{\rm SO}}

\newcommand{\Def}{{\bf Definition}}
\newcommand{\Rem}{{\bf Remark}}
\newcommand{\Prop}{{\bf Proposition}}
\newcommand{\Le}{{\bf Lemma}}
\newcommand{\Th}{{\bf Theorem}}

\newcommand{\Cor}{{\bf Corollary}}

\newcommand{\Exa}{{\bf Example} }

\newcommand{\s}{{\mathfrak {sl}}}

\newcommand{\ot}{{\mathfrak o}}
\newcommand{\syp}{{\mathfrak {sp}}}

\newcommand{\p}{{\mathfrak p}}

\newcommand{\ka}{{\mathfrak k}}

\newcommand{\g}{{\mathfrak g}}
\newcommand{\h}{{\mathfrak h}}

\newcommand{\fa}{\quad {\rm for \,\, all}\,\,\,}
\newcommand{\with}{\quad {\rm with}\,\,\,}
\newcommand{\Ima}{\quad {\rm Im}\,\,\,}
\newcommand{\sig}{\,\, {\rm sig}\,\,}

\begin{document}
\thispagestyle{empty}
\phantom{}
\vspace{1 cm}
{\bf  \Large Round about Theta. Part I Prehistory}\\

Rolf Berndt\\

There is a huge amount of work on different kinds of theta functions, the theta correspondence, cohomology classes 
coming from special Schwartz classes via theta distribution, and much more. The aim of this text is to try to find joint 
construction principles while often leaving aside relevant but cumbersome details.\\

The next steps after this prehistoric Part I will be directed to a description of the Howe operators introduced by 
Kudla and Millson and their special Schwartz forms and classes. This has as attractor the fact that the 
modular and automorphic forms arising naturally in context with these classes find very nice geometric interpretations 
of their Fourier coefficients and thus lead to an intriguing intertwining of elements of representation theory with 
algebraic and arithmetic geometry.\\

The presentation here is in the spirit of my book on representations of linear groups. 
Though it may be seen as just another chapter, it has it has its own raison d'\^etre and can be read independently.\\

\tableofcontents

\section {Riemann and Jacobi Theta Series}

\subsection {Weil Representation} 

{\bf 1.1.1} Our fundamental object is the symplectic group
$$
\hat G = \Spe(n,\R) := \{ g \in \M(2n,\R); {}^tgJg = J := \begin{pmatrix} &1_n\\-1_n& \end{pmatrix} \},
$$
i.e.~the group with elements
$$
g = \begin{pmatrix} A&B\\C&D \end{pmatrix}, A,B,C,D \in \M(n,\R), {}^tAD -{}^tCB = E_n, {}^tAC = {}^tCA, {}^tBD = {}^tDB,
$$
and its projective (Segal-Shale-)Weil or oscillator \re $\omega $ given as its Schr\"odinger model on the space $\Hi = L^2(\R^n)$ 
by the prescription
$$
\begin{array}{rcl}
\omega (d(A))f(x) & = &\mid \det A \mid ^{1/2} f({}^t\!Ax) \fa d(A) :=  \begin{pmatrix} A&\\&{}^tA^{-1} \end{pmatrix}, A \in \GL(n,\R), \\
\omega (n(B))f(x) & = &e^{\pi i{}^txBx}f(x) \fa n(B) :=  \begin{pmatrix} 1&B\\&1 \end{pmatrix}, B \in \Sym(n,\R),\\
\omega (J)f(x) & = &\gamma \hat f(x), \hat f(x) := \int_{\R^n}f(y)\,e^{2\pi i\ta yx}dy,
\end{array}
$$
where $\gamma $ will be specified later.
This projective \re corresponds to a \re $\tilde \omega $ of the twofold cover of $\hat G$, the metaplectic group $\tilde {\hat G} = \Mp(n,\R)$ 
with elements $(g,t), \, g \in \hat G, t^2 = s(g)^{-1},$ where 
$s(g)$ as specified in [LV] p.70 will not be needed at the moment. \\

{\bf 1.1.2} The Lie algebra of the symplectic group is 
$$
\hat \g = \syp(n,\R) := \{ X = \begin{pmatrix} A&B\\C&-{}^t\!A \end{pmatrix}; A,B,C \in \M(n,\R), B = {}^tB, C = {}^tC \}.
$$
$\syp$ has dimension $2n^2+n$ and Cartan decomposition
$$
\syp = \ka + \p
$$
with
$$
 \ka := \{ \begin{pmatrix} A&B\\-B&A \end{pmatrix}; {}^tB = B, {}^t\!A = - A \}, \,\,
 \p := \{ \begin{pmatrix} A&B\\B&-A \end{pmatrix}; {}^tB = B, {}^t\!A =  A \}.
$$
The complexification $\g_c$ of $\g = \syp$ has the $\Ad J-$eigenspace decomposition
$$
\g_c = \syp^{(1,1)} + \syp^{(2,0)} + \syp^{(0,2)}
$$
with
$$
\begin{array}{rcl}
\syp^{(1,1)} & := &  \{ \begin{pmatrix} A&B\\-B&A \end{pmatrix}; A,B \in \M(n,\C), {}^tB = B, {}^t\!A = - A \}, \\
\syp^{(2,0)} & := &  \{ \begin{pmatrix} A&iA\\iA&-A \end{pmatrix};  A \in \M(n,\C), {}^t\!A =  A \}, \\
\syp^{(0,2)} & := &  \{ \begin{pmatrix} A&-iA\\-iA&-A \end{pmatrix};  A \in \M(n,\C), {}^t\!A =  A \}. 
\end{array}
$$
We take over the notation from Adams ([Ad] p.466)
$$
A_{ij} := \begin{pmatrix} E_{ij}&\\&-E_{ji} \end{pmatrix}, \,\,U^+_{ij} := \begin{pmatrix} 0&B_{ij}\\0&0 \end{pmatrix},\,\,
 U^-_{ij} := \begin{pmatrix} 0&0\\B_{ij}&0 \end{pmatrix}
$$
where for $1\leq i,j \leq n$ $E_{ij}$ is the elementary matrix with zero entries except there is 1 in the $i$th row and $j$th column  
and $B_{ij} := E_{ij}+E_{ji}$ for $i \not= j$, $B_{ii} := E_{ii}.$\\

{\bf 1.1.3} For $X \in \g_c$ we denote by $\hat X$ its operator in the derived \re of the Weil \re $\omega, $ i.e.~we put
$$
\hat X f(x) := \frac{d}{dt}\vert_{t=0} (\omega (\exp tX)f)(x) \fa f \in \Ss(\R^n).
$$
One easily comes to
$$
\hat A_{jk} = x_j \partial _k + \delta _{jk}/2, \,\,\,\hat U^+_{jk} = 2\pi i x_jx_k, j\not=k,\,\, \hat U^+_{jj} = \pi i x^2_j,
$$
and, using appropriate commutation formulae like $[U^+_{jj},U^-_{jj}] = A_{jj},$ 
$$
\hat U^-_{jj} = - 1/(4\pi i) \partial ^2_j,\,\,\hat U^-_{jk} = - 1/(2\pi i) \partial _j\partial _k \fa j \not=k.
$$
Hence, for the complex algebras $\syp^{(1,1)} = \ka_c, \,\,\syp^{(2,0)} =:\p^+,\,\syp^{(0,2)} =: \p^-,$ respectively generated by 
$$
A_{jk}-A_{kj},\,\,U^+_{jk}-U^-_{jk}
$$
and
$$
 (1/2)(\mp i(A_{jk}+A_{kj})/(1+\delta _{jk}) + U^+_{jk} + U^-_{jk}) =: \check U^\pm_{jk},
$$
one has
$$
\begin{array}{rcl}
\hat { \check U}^\pm_{jj} & = & (1/2)(\mp i(x_j\partial _j + (1/2)) + \pi ix_j^2 - (1/4\pi i)\partial ^2_j),\\
\hat { \check U}^\pm_{jk} & = & (1/2)(\mp i(x_j\partial _k + x_k\partial _j) + 2\pi ix_jx_k - (1/2\pi i)\partial _j\partial _k).
\end{array}
$$
In $\p$ one has the Cartan algebra $\h := \,<A_{jj}>_{j=1,\dots,n}$ and (among others) the relations
$$
[A_{jj},U^\pm_{jj}] = \pm2U^\pm_{jj},\,\,\,[U^+_{jj},U^-_{jj}] = A_{jj}.
$$
We use the Cayley transformation, i.e.~conjugation by
$$
c = (1/\sqrt 2) \begin{pmatrix} 1_n&i_n\\i_n&1_n \end{pmatrix},
$$
to introduce
$$
H_j := cA_{jj}c^{-1} = -i(U^+_{jj} - U^-_{jj})
$$
which obeys the relation
$$
[H_j,\check U^\pm_{jj}] = \pm 2\check U^\pm_{jj}.
$$
We get
$$
\hat H_j = \pi x^2_j - (1/4\pi)\partial _j^2.
$$
{\bf 1.1.4} Now we can see that for the Gaussian
$$
\varphi _0 (x) := e^{-\pi \Sigma x^2_j} 
$$
one has
$$
\begin{array}{rcl}
\hat H_j \varphi _0 & = & (1/2) \varphi _0,\\
\hat {\check U}^+_{jj} \varphi _0 & = & (i/2)(4\pi x^2_j - 1) \varphi _0,\\
\hat {\check U}^+_{jk} \varphi _0 & = & 4\pi i x_jx_k \varphi _0, \fa j\not=k,\\
\hat {\check U}^-_{jk} \varphi _0 & = & 0, \fa j\not=k,\\
\hat {\check U}^-_{jj} \varphi _0 & = & 0.
\end{array}
$$
This shows that $\varphi _0$ is annihilated by all elements of $\syp^{(0,2)}$ and reproduced with eigenvalue 1/2 by all elements 
of $\syp^{(1,1)},$ i.e.~$\varphi _0$ is a {\it vacuum vector} for $\omega :$ a vector of lowest weight 1/2 for the Weil \re $\omega $ 
in its Schr\"odinger model. \\

\subsection {Riemann Thetas}

{\bf 1.2.1}  There is a standard way to construct a modular form which in this case comes out like this: 
One applies the Weil \re $\omega $ to $\varphi _0$ and averages over all $\ell \in \Z^n$ to get a function 
$$
\Phi _\theta (g) := \sum_{\ell \in \Z^n} (\omega (g)\varphi _0)(\ell ),\,\,g \in \hat G,
$$
which can be proven to be invariant under the theta subgroup $\Gamma _\theta $ of $\Spe(n,\R).$ 
And up to an automorphic factor this function can be identified with $\theta,$ the (zero value) of the Jacobi theta function: 
With some more (but not all) details, this means the following.\\

{\bf 1.2.2} We have the transitive action of $\Spe(n,\R)$ on the Siegel half space
$$
\Ha_n := \{\hat \tau \in \M(n,\C); \,\,{}^t\hat \tau  = \hat \tau , \Ima  {\hat \tau} > 0\}
$$
given by
$$
 g(\hat \tau ) := (A\hat \tau +B)(C\hat \tau +D)^{-1} \fa g = \begin{pmatrix} A&B\\C&D \end{pmatrix}.
$$ 
We take an element $g = g_{\hat \tau} \in \hat G$ such that 
$$
g_{\hat \tau} (i_n) = \hat \tau =: \hat u + i\hat v,
$$
namely, using the notation introduced above
$$
g_{\hat \tau } = n(\hat u)d(A), \with A{}^t\!A = \hat v.
$$
Then we get
$$
(\omega (g_{\hat \tau })\varphi )(x) = \vert \det \hat v \,\vert^{1/4} e^{\pi i {}^tx \hat \tau x}
$$
and hence
$$
\Phi _\theta (g_{\hat \tau} ) :=  \vert \det \hat v \,\vert^{1/4} \sum_{\ell \in \Z^n} e^{\pi i {}^t\ell \hat\tau  \ell}.
$$
Here we find the standard theta series
$$
\theta (\hat \tau) :=  \sum_{\ell \in \Z^n} e^{\pi i {}^t\!\ell \hat\tau  \ell}
$$
{\bf 1.2.3} The fact that $\varphi _0$ is a lowest weight vector annihilated by $\syp^{(0,2)}$ tanslates into the fact that 
$\theta $ is a holomorphic function in $\hat \tau \in \Ha_n$ and the fact (which is not so easy to prove) 
that $\Phi _\theta$ as a function of $g \in \hat G$ is invariant under the theta group $\Gamma _\theta $ translates into the 
automorphic functional equation
$$
\theta (\gamma \hat \tau ) = \varepsilon (\gamma ) \det(C\hat \tau + D)^{1/2}\theta (\hat \tau )
$$
for $ \gamma \in \Gamma _\theta,$ the group of elements 
$$
 \gamma =  \begin{pmatrix} A&B\\C&D \end{pmatrix} \in \Gamma := \Spe(n,\Z) ,
$$
where ${}^tCA$ and ${}^tBD$ both have even diagonal entries. $\varepsilon (\gamma ) $ is a character of $\Gamma _\theta $ as defined in 
[LV] p.166. In particular for the group 
$$
\Gamma ^0_0(2) := \{\gamma  = \begin{pmatrix} A&B\\C&D \end{pmatrix} \in \Gamma ;\,\,\, B \equiv C \equiv  0 \mod 2 \}
$$
one has
$$
\epsilon (\gamma )^2 = (\frac{-1}{\det D}).
$$ 
This statement is Theorem 2.2.37 in [LV]. We shall call this $\theta $ as a function on the Siegel half space the Riemann theta function, 
though this name is used by Mumford also for the more general function which we introduce now and then call Jacobi theta function. \\

\subsection {Jacobi Thetas}

{\bf 1.3.1} For $\tau \in \Ha_n$ and $z \in \C^n$ we get the Jacobi theta function
$$
\theta (\tau ,z) :=  \sum_{\ell \in \Z^n} e^{\pi i ({}^t\!\ell \hat\tau  \ell + 2{}^tz\ell)}.
$$
{\bf 1.3.2} Here we have to extend the symplectic group $\hat G$ to its semidirect product with an appropriate Heisenberg 
group $\Heis (\R^n)$ to come to the Jacobi group $\hat G^J.$ As a set one has $\Heis (\R^n) = \R^{2n+1}$ and all multiplication laws 
are fixed by the embedding into the symplectic group $\Spe(n+1,\R)$ given by
$$
\Heis (\R^n) \ni (\lambda ,\mu , \kappa ) \longmapsto 
\begin{pmatrix} 1_n&&&\mu \\{}^t\lambda &1&{}^t\mu &\kappa \\&&1_n&-\lambda\\&&&1  \end{pmatrix} ,
$$
$$
\Spe(n,\R) \ni M =  \begin{pmatrix} A&B\\C&D \end{pmatrix} \longmapsto 
\begin{pmatrix} A&&B&\\&1&&\\C&&D&\\&&&1 \end{pmatrix} .
$$
We write
$$
g = (p,q,\kappa )M \,\,\,{\rm or}\,\,\,g = M(\lambda ,\mu ,\kappa ) \,\in G^J(\R^n).
$$
$\Heis(\R^n)$ acts on $\R^{2n}$ via $(x,y)\longmapsto (x+\lambda ,y+\mu )$ and
$G^J$ acts on $\Ha_n\times \C^n$ via
$$
(\tau ,z)\longmapsto g(\tau ,z) := (M(\tau ),(z + \tau \lambda + \mu )(C\tau + D)^{-1})
$$
where $g = M(\lambda ,\mu ,\kappa ) \in G^J(\R^n), \tau \in \Ha_n, z \in \C^n$. For $g = (p,q,\kappa )M$ one has
$$
g(i_n,0) = (\tau = M(i_n), z = \tau p + q).
$$
{\bf 1.3.3} The construction of the Weil \re usually goes via the
standard \re of the Heisenberg group which is the Schr\"odinger \re in the space $\Hi = L^2(\R^n)$ 
for real non-zero $m$ and $(\lambda ,\mu ,\kappa ) \in \Heis(\R^n)$ given by 
$$
(\pi_S^m(\lambda ,\mu ,\kappa )f)(x) := e^m(\kappa  + (2{}^tx+{}^t\lambda)\mu )f(x+\lambda ) \,\,\fa f \in \Hi.
$$
Then one has the Schr\"odinger-Weil \re $\pi_{SW}$ of $G^J$ given by
$$
\pi^m_{SW}((p,q,\kappa )M) := \pi_S^m((p,q,\kappa ) \omega (M)
$$
{\bf 1.3.4} It is not difficult to verify that the vacuum vector of the Weil \re $\varphi _0^m(x) = e^{\pi i m{}^t xx}$ is also a vacuum vector 
of the Schr\"odinger-Weil \re and one can use it again as done above: For $M_{\hat \tau} = n(\hat u)d(A)$ we get
$$
(\pi_{SW}((p,q,\kappa )M_{\hat \tau }) \varphi _0^m)(x) = 
\vert \det \hat v\,\vert^{1/4}e^m(\kappa + {}^tp\hat \tau p + {}^tpq) e^{\pi i m({}^tx\hat \tau x + 2({}^tp\hat \tau + {}^tq)x)} 
$$
and
$$
\Phi _\theta ((p,q,\kappa )M_{\hat \tau }) = 
\sum_{\ell \in \Z^n} \vert \det \hat v\,\vert^{1/4}e^m(\kappa + {}^tp\hat \tau p + {}^tpq) e^{m\pi i ({}^t\ell\hat \tau \ell + 2({}^tp\hat \tau + {}^tq)\ell)} .
$$
With $z = p\hat \tau + q,$ for $m = 1/2,$ up to a factor we find the Jacobi theta function
$$
\theta (\tau ,z) :=  \sum_{\ell \in \Z^n} e^{\pi i ({}^t\!\ell \hat\tau  \ell + 2{}^tz\ell)}.
$$
The properties and the functional equation of this function and its generalizations are discussed 
with the appropriate details in the books by Igusa ([Ig] p.48f) and by Mumford ([MuIII] p.142). Here 
we only will record the following observation. \\

{\bf 1.3.5} \Rem: We introduced the Jacobi group $G^J(\R^n)$ 
as a subgroup of the symplectic group $G^\ast  := \Spe(n+1;\R)$. We have 
$$
g^\ast (i_{n+1}) = (A^\ast i_{n+1} + B^\ast )(C^\ast i_{n+1} + D^\ast )^{-1} \,\,\,{\rm for}\,\,
g^\ast = \begin{pmatrix} A^\ast &B^\ast \\C^\ast &D^\ast \end{pmatrix}.
$$
If we specialize this for $g^\ast = (p,q,\kappa )M_{\hat \tau },$ we get
$$
((p,q,\kappa )M_{\hat \tau }) (i_{n+1}) =  \begin{pmatrix} \hat \tau &p\hat \tau +q \\{}^tp \hat \tau +{}^tq&{}^tp \hat \tau p + {}^tqp + \kappa + i \end{pmatrix}.
$$
And if we specialize the standard theta series for $G^\ast $
$$
\theta_{n+1} (\tau ^\ast ) = \sum_{\ell^\ast \in \Z^{n+1}} e^{\pi i {}^t\!\ell^\ast \tau ^\ast \ell ^\ast }
$$
for 
$$
 \tau^\ast _0 = \begin{pmatrix} \hat \tau &z \\{}^tz&{}a\end{pmatrix}, \,\, a = \ta p \hat \tau p + \ta qp + \kappa + i 
$$
with ${}^t\!\ell^\ast = ({}^t\!\ell,l), \,\ell  \in \Z^n, \,l \in \Z,$ we get
$$
\theta_{n+1} (\tau ^\ast_0 ) = \sum_{l\in\Z} e^{\pi ia l^2}\sum_{\ell\in\Z^n} e^{\pi i ({}^t\!\ell \hat \tau \ell + 2l\ta z\ell)}.
$$
Up to the factor $l \in \Z$ in the exponent, we find again the Jacobi theta series. This fits into the framework of 
the Fourier-Jacobi expansion of a Siegel modular form which is introduced (for the lowest dimensional case) in 
[EZ] p.72f. From here we can easily take over that each coefficient 
$$
\phi _{l^2}(\tau ,z) := \sum_{\ell\in\Z^n} e^{\pi i ({}^t\!\ell \hat \tau \ell + 2l\ta z\ell)}
$$
has the transformation property of a Jacobi form and using the operator $U_l$ defined in [EZ] p.41
(multiplication of the $z-$variable by $l$) we can even write
$$
\phi _{l^2}(\tau ,z) = U_l \,\theta (\tau ,z).
$$
{\bf 1.3.6} There are many ways to introduce more general functions of this type. As one can well imagine, all this 
generalizes rather easily if one takes a rational symmetric $h\times h-$matrix $S$ belonging to a positive 
definite quadratic form.  We follow [MuIII] p.96f: \\

{\bf 1.3.7} \Def: Let $S \in \Sym_h(\Q)$ be positive definite, $T \in \Ha_n$ and $Z \in M_{n,h}(\C).$ Then we put
$$
\theta ^S(T,Z) := \sum_{N \in M_{n,h}(\Z)} e^{\pi i \Tr ({}^tNTNS + 2{}^tNZ)}.
$$
As it is rather easy to see that for $M,N \in \M_{g,h}(\Z)$ one has
$$
\theta ^S(T,Z + TMS + N) e^{\pi i \Tr({}^tMTMS + 2{}^tMZ)} = \theta ^S(T,Z)
$$
one is lead to suggest that $\theta ^S$ is just a Jacobi theta series for a more general situation, 
namely for the complex torus
$$
M_{n,h}(\C)/(T\M_{n,h}(\Z)S + \M_{n,h}(\Z)).
$$
To see this we use the identifications given by
$$
\M_{n,h}(\C) \longrightarrow \C^{nh},\,\,\, Z = (Z_1,\dots,Z_h) \longmapsto z = \begin{pmatrix} Z_1 \\ \cdot \\ \cdot \\ Z_h \end{pmatrix}, 
$$
and
$$
\tau := T \otimes S = \begin{pmatrix} TS_{11} & \cdot & \cdot & TS_{1h} \\ \cdot &&&\cdot \\\\ \cdot &&&\cdot \\TS_{h1}&\cdot &\cdot &TS_{hh}\end{pmatrix}.
\,\,\in \M_{nh,nh}(\C).
$$ 
One has to check that $W = TZS$ translates into $w = \tau z$ and that one has $\Tr {}^tWZ = {}^twz$. 
Then we get (Lemma 6.2 in [MuIII])
$$
M_{n,h}(\C)/(T\M_{n,h}(\Z)S + \M_{n,h}(\Z)) \cong \C^{nh}/(\tau \Z^{nh} + \Z^{nh}).
$$
Hence we can see that we have
$$
\theta ^S(T,Z) = \sum_{n\in\Z^{nh}} e^{\pi i({}^tn\tau n + 2{}^tnz)} = \theta_{nh} (\tau ,z).
$$
Moreover one can see without too much trouble ([MuIII] Corollary 6.6):\\
If $S = \begin{pmatrix} d_1& & &0 \\ &d_2& & \\ &&\cdot& \\ 0&&&d_h \end{pmatrix}$ and $Z = (z_1,\dots,z_h),$ then
we have
$$
\theta ^S(T,Z)  = \Pi _{i=1}^h \theta_n (d_iT,z_i).
$$
\vspace{.1cm}

{\bf 1.3.8} Thetas belonging to (positive definite) quadratic forms are still more widely generalized by considering 
spherical harmonic polynomials as coefficients of the exponentials in the series. This is treated for instance in 
[MuIII] p.145ff. Here we shall  come to this later.\\ 

\section {Hecke and Siegel Theta Series}

\subsection {Hecke Thetas}

{\bf 2.1.1} It is immediate that one has a convergence problem if one tries to consider thetas for quadratic forms which are 
not positive definite. For instance, for $\tau \in \Ha := \Ha_1$ the sum
$$
\Sigma _{x_1,x_2 \in \Z}e^{2\pi i \tau (x_1^2 - 12x_2^2)}
$$
has no sense. It was Hecke on his way to associate modular forms to real quadratic fields (and a bit later Schoeneberg) 
who achieved substantial progress in this topic:\\
Let $K = \Q(\sqrt D)$ be a real quadratic number field with diiscriminant $D$ and $\ot := \ot_K$ as its ring of integers. 
For $Q \in \N$ and $\alpha \in \ot$ Hecke defines in [H1] and [H2] the functions of $\tau \in \Ha$
$$
\vartheta (\tau ; \alpha , Q\sqrt D) := \sum_{(\mu )} \sgn \mu \,\,e^{2\pi i\tau \frac{\vert \mu \mu '\vert}{QD} }
$$
and
$$
\vartheta_+ (\tau ; \alpha , Q\sqrt D) := \sum_{(\mu ), \mu \mu '>0} \sgn \mu \,\,e^{2\pi i\tau \frac{\vert \mu \mu '\vert}{QD} }
$$
Here the prime ' indicates the conjugate in $K$ and the summation $\sum_{(\mu )}$ is meant over a family of elements $\mu \in \ot,$ 
which are congruent $\mod Q\sqrt D$ to $\alpha $ and not associated, i.e.~do not differ by a unit $\mod Q\sqrt D$ as a factor. 
Hecke's main theorem in this context is a transformation formula $\tau \longmapsto -1/\tau$ for $\vartheta _+,$ 
which is the essential to show that $\vartheta _+$ is a modular form. Without going into further details we state 
that the idea for his proof is to use the already known transformation property of a standard theta function in two variables. \\

{\bf 2.1.2} Hecke discusses the example $\vartheta _+(\tau ; 1, \sqrt {12})$. He shows that one has
$$
\begin{array}{rcl}
\vartheta _+(\tau + 1; 1, \sqrt {12})  & = & e^{2\pi i/12} \vartheta _+(\tau ; 1, \sqrt {12}),\\
\vartheta _+(-1/\tau ; 1, \sqrt {12})  & = & - i\tau  \vartheta _+(\tau ; 1, \sqrt {12}), 
\end{array}
$$
and one has the nice relation to the Delta function
$$
 \vartheta _+(\tau ; 1, \sqrt {12}) = (\Delta (\tau ))^{1/12} = e^{2\pi i\tau /12} \Pi _{i=1}^\infty (1-  q^n)^2, \,\,q := e^{2\pi i\tau }.
$$
\subsection {Siegel Thetas}

{\bf 2.2.1} As we saw, Hecke solved the convergence problem for the theta series for indefinite quadratic forms by 
summing only over those elements such that the form has positive values. Siegel had the idea to use the {\it majorant} 
of a quadratic form to associate to the form a convergent series for which (in [S1] and [S2]) he also could prove a modular property.\\ 

{\bf 2.2.2} One starts with the quadratic form belonging to a non-degenerate symmetric matrix $S \in \M_n(\R)$ with signature $\sig S = (p,q)$
$$
S[x] := \ta xSx, 
$$ 
where as above $x$ is a column. We know that one can find a matrix $C \in \GL(n,\R)$ such that
$$
S[C] = \ta C S C =  \begin{pmatrix} 1_p&  \\ &-1_q \end{pmatrix} =:S_0,
$$
i.e.~with $x = Cy$ one has
$$
S[Cy] = S_0[y] = y_1^2 +\dots+ y^2_p - (y^2_{p+1} +\dots+y^2_{p+q}) =: {y'}^2 - {y''}^2.
$$
Siegel now uses the notion of the majorant of $S[x]$ which goes back to Hermite and is a positive definite quadratic form, 
say $P[x]$, such that $P[x] \geq S[x]$ for all $x \in \R^n$. With $C$ as above, we take $P := (C\ta C)^{-1}$ and get
$$
P[x] = \ta x(C\ta C)^{-1}x = \ta yy.
$$
In [S2] 1. Siegel shows that $P$ belongs to a majorant of $S[x]$ if and only if $P$ fulfills the two conditions 
$$
PS^{-1}P = S, \,\, \ta P = P >0.
$$
Moreover, Siegel parametrizes the set $\Pp := \Pp(S)$ of these matrices $P$ and shows that the orthogonal group 
$$
\Or := \Or(S) = \{ A \in \M_n(\R); \,\, \ta ASA = S \}
$$
via $(A,P) \longmapsto P[A]$ acts transitively on $\Pp$. We will come back to this later but now can give 
Siegel's definition of his Theta function:\\

{\bf 2.2.3} \Def: One takes $\tau = u + iv \in \Ha,$ $R := uS + ivP$ and puts
$$
\theta (\tau ) := \theta (\tau ,P) := \sum_{x \in \Z^n} e^{2\pi i R[x]}.
$$

This definition makes sense because$\Ima R = vP$ is positive definite.  For $a \in \Q; as \in \Z$ where $s := \det S$ Siegel also looks 
at the variant
$$
\theta_a (\tau ) := \theta_a (\tau ,P) := \sum_{x \in \Z^n} e^{2\pi i R[x +a]}.
$$
{\bf 2.2.4} This function $\theta_a $ is not a holomorphic function in $\tau $ but has 
a modular behaviour with respect to certain modular substitutions
$\tau \longmapsto \hat \tau = (a\tau + b)(c\tau + d)^{-1}$ which is given in Hilfssatz 1 in [S2]. 
The proof again uses the Poisson summation formula.
We will not repeat this here but just indicate that an automorphic factor of type
$$
(c\tau + d)^{-p/2}(c\bar \tau + d)^{-q/2}
$$
comes in. The dependence on $P$ resp.~on appropriate parameters for $P$ will be discussed later.\\

{\bf 2.2.5} There are several generalizations of Siegel's definition, see, for instance, Vign\'eras [Vi] and in particular 
Borcherds [Bo]. Moreover, there is an extension in the direction of Jacobi thetas by O. Richter [Ri]: \\

{\bf 2.2.6} \Def: Let $S \in \Sym_m(\Z)$ be an invertible matrix with even diagonal entries 
with $\sig S = (p,q)$ and such that $qS^{-1}$ for $q \in \N$ is integral and even. 
Let $P$ be a majorant of $S, \tau = u + iv \in \Ha_n, \zeta \in \M_{m,j}(\Z),$ and $Z  \in \M_{j,n}[\C).$ Then one puts
$$
\theta _{S,P,\zeta }(\tau ,Z) := \sum_{N \in \M_{m,n}} e^{\pi i \Tr(S[N] u + i P[N] v + 2\ta NS\zeta Z)}.
$$
For $\zeta $ such that $S\zeta = P\zeta ,$\, Richter proves a transformation formula concerning
$$
\Gamma ^{(n)}_0 := \{ \begin{pmatrix} A&B \\ C&D \end{pmatrix} \in \Spe(n,\Z); \,\, C \equiv 0 \mod q \}.
$$
The definition of these theta series is quite natural in the context the autors had. But they also have 
a \re theoretic background as we will try to elucidate in the sequel.\\

\section {A Dual Pair and Siegel Thetas}

\subsection {Dual Pairs}

It was Howe who (in [Ho1]) coined the following notion and later contributed essential parts 
of its discussion.\\

{\bf 3.1.1} \Def: A dual reductive pair is a pair of subgroups $(G,G')$ in a symplectic group $\hat G =\Spe(n,\R)$ such that\\
i) $G$ is the centralizer of $G'$ in $\hat G$ and $G'$ is the centralizer of $G$ in $\hat G.$\\
ii) The actions of $G$ and $G'$ on $\hat V := \R^{2n}$ are completely reducible (i.e.~every invariant subspace has an invariant complement).\\

This is only a special case: here one can also replace $\R$ by more general fields. 
One defines irreducible pairs as those where one can not decompose $\hat V$ as the direct sum 
of two symplectic subspaces each of which is invariant under both $G$ and $G'.$ 
There is the classification of irreducible pairs done in [MVW]. We won't go into this but just 
point out that these pairs provide the background for a lot of important 
relations between different kinds of automorphic forms. Roughly, this goes like this: 
The Weil \re $\omega $ of $\hat G$ restricts to \res of the subgroups $G$ and $G'$ and to $G\times G'$. 
If one has a decomposition of $\omega $ where to an irreducible  \re of $G$ corresponds 
exactly one  irreducible  \re of $G'$,  one can hope for a correspondence between automorphic forms belonging 
to these \res. With more details one has the Howe conjecture making precise statements in this 
direction. For the moment we will use a small part of the picture to make reappear the Siegel thetas.\\

{\bf 3.1.2} We take the orthogonal group $G = \Or(p,q) \cong \Or(S)$ belonging to the non-degenerate symmetric matrix 
$S \in \M_n(\R)$ with signature $\sig S = (p,q).$ Then one can verify easily that $G$ together with 
$G' = \SL(2,\R)$ is a dual pair in $\hat G = \Spe(n,\R).$\\ 

{\bf 3.1.3} \Rem: We use the embeddings
$$
G = \Or(p,q) \ni A \longmapsto \hat A := \begin{pmatrix} \ta A^{-1}& \\ &A \end{pmatrix} \in \hat G = \Spe(n,\R)
$$  
and with $S_0 = \begin{pmatrix} 1_p& \\ &-1_q \end{pmatrix}$ 
$$  
G' = \SL(2,\R) \ni M := \begin{pmatrix} a&b \\c&d \end{pmatrix} \longmapsto \hat M :=
 \begin{pmatrix} a1_n&bS_0 \\cS_0^{-1} &d1_n \end{pmatrix}  \in  \hat G = \Spe(n,\R).
$$ 
These embeddings come as special cases from the following more general consideration: \\

{\bf 3.1.4} For the symplectic space $V' \simeq \R^{2m}$ with the action of $G' = \Spe(m,\R)$ we take as basis $e_1,\dots,e_m,e'_1,\dots,e'_m$ 
such that for all $j=1,\dots,m$ and $J_m = \begin{pmatrix} &1_m \\-1_m & \end{pmatrix}$ one has $J_me_j = e'_j$ and $J_me'_j = -e_j.$ 
For the orthogonal space $V \simeq \R^n$ with the action of $G = \Or(p,q), p+q =n,$ we take as basis $v_1,\dots,v_p,v_{p+1},\dots,v_{p+q}$ 
such that for all $\alpha = 1,\dots,p$ and $\nu = p+1,\dots,p+q$ and $S_0 = \begin{pmatrix} 1_p&\\ &-1_q \end{pmatrix}$ one has
$S_0v_\alpha = v_\alpha $ and $S_0v_\nu = - v_\nu .$ Hence $\hat V := V \otimes V' \simeq \R^{2mn}$ is a symplectic space 
with basis $e_j\otimes v_\alpha ,e_j\otimes v_\nu $ and $e'_j\otimes v_\alpha , - e'_j\otimes v_\nu,$ i.e.~ 
$\hat e_1,\dots,\hat e_{mn},\hat e'_1,\dots,\hat e'_{mn}$ where
$$
\hat e_1 := e_1\otimes v_1,\dots,\hat e_m := e_m\otimes v_1,\dots,\dots,\hat e_{mn} := e_m\otimes v_{p+q}
$$
and
$$
\hat e'_1 := e'_1\otimes v_1,\dots,\hat e'_{mp} := e'_m\otimes v_q,\hat e'_{mp+1} := - e'_1\otimes v_{p+1},
\dots,\hat e'_{mn} :=  - e'_m\otimes v_{p+q}.
$$
In particular for $m = 1$ we have $\hat V $ with basis
$$
\hat e_j:= e_1\otimes v_j,\,\,\hat e'_j:= e'_1\otimes S_0v_j,\fa j = 1,\dots,n.
$$
The action of $G$ on $V$ and of $G'$ on $V'$ induce naturally actions on $\hat V,$ in particular those described in the remark above.\\ 

\subsection {Siegel Thetas as Special Values of Riemann Thetas}

{\bf 3.2.1} We have the standard compact subgroups
$$
K := \Or(p) \times \Or(q) ,\,\, K' := \SO(2),
$$ 
and
$$
\hat K := \{ \begin{pmatrix} A&B \\-B &A \end{pmatrix};\,\, \ta AA + \ta BB = 1_n, \,\ta AB = \ta BA \,\} \simeq \Un(n).  
$$
There are the standard maps to the associated homogeneous spaces
$$
\hat G \longrightarrow \hat G/\hat K = \Ha_n; \,\, \hat g =  \begin{pmatrix} A&B \\C&D \end{pmatrix}  \longmapsto \hat g(i_n) =:\hat \tau = \hat u + i\hat v,
$$
$$
G' \longrightarrow G'/K' = \Ha;\,\,g =  \begin{pmatrix} a&b \\c&d \end{pmatrix} \longmapsto g(i) = \tau = u+iv,
$$
and (without a big loss of generality) restricting to $S = S_0$
$$
G \longrightarrow G/K =: \D ;\,\, A \longmapsto (A\ta A)^{-1} =: P.
$$
This homogeneous space has different realizations which we will discuss later. Here we refer to our remarks in 2.2.2 where, 
following Siegel, we introduced $\D = \Pp$ as the set of majorants of $S =S_0$. The embedding 
$G \times G' \longrightarrow \hat G$ induces a map
$$
\D \times \Ha \longrightarrow \Ha_n; \,\, (P,\tau ) \longmapsto uS_0 + ivP =:\hat {\tau} _{P,\tau } =:\hat {\tau _0},
$$
which is a consequence of 
$$
\hat g_\tau(i_n) = uS_0 + vi_n \,\,{\rm and}\,\,\hat A(\hat \tau ) = \ta A^{-1}\hat \tau A^{-1}.
$$ 
{\bf 3.2.2} \Rem: If we specialize the variable $\hat \tau $ in the standard Riemann theta series $\theta (\hat \tau )$ for $\Spe(n,\R)$ to 
$\hat \tau = \hat \tau _0,$ 
we recover the Siegel theta series
$$
\vartheta (\tau ,P) = \sum_{\ell \in \Z^n} e^{\pi i \ta \ell (uS_0 + ivP)\ell}.
$$
In a parallel way, one can take the vacuum vector $\varphi _0(x) = e^{\pi \Sigma x_j^2}$ for the Weil \re $\omega $ in the Schr\"odinger model 
and apply the restriction of $\omega $ to $G\times G'$ to construct a function on $G\times G'$ with certain invariance properties. 
This way we come to
$$
\omega (\hat A \cdot \hat {g_\tau })\varphi _0(x) = v^{n/4}e^{\pi i \ta x(u S_0 + iv \ta A^{-1} A^{-1})x}.
$$

\subsection {Siegel Theta and its Representation}

In [S2] Siegel  uses these theta series to study the diophantine problem of integral solutions $x \in \Z^n$ of the 
quadratic equation 
$$
S[x+a] = t.
$$
Here we won't go into this interesting topic but analyze a bit the relation of the Siegel theta series to the \re theory of the 
two groups $G$ and $G'$ going into our construction.\\

{\bf 3.3.1} In 1.1.2 we discussed the Lie algebra $\hat \g$ of the symplectic group $\hat G = \Spe(n,\R).$ As a special case we 
have $\g' = \Lie G', G' = \SL(2,\R)$ with
$$
\g' = < F := \begin{pmatrix} &1 \\& \end{pmatrix}, \,\,G := \begin{pmatrix} & \\1& \end{pmatrix}, \,\,H := \begin{pmatrix} 1& \\&-1 \end{pmatrix} >
$$ 
and the relations
$$
[H,F] = 2F,\,\,[H,G] = - 2G,\,\, [F,G] = H.
$$
The complexification is given by
$$
\g'_c = <Z := -i(F - G) = \begin{pmatrix} &-i \\i& \end{pmatrix}, X_{\pm} := (1/2)(H \pm i(F+G)) = 1/2\begin{pmatrix} 1&\pm i \\\pm i&-1 \end{pmatrix} > 
$$
and the relations
$$
[Z,X_{\pm}] = \pm 2X_{\pm},\,\,[X_+,X_-] = Z.
$$
{\bf 3.3.2} If we use the embedding of $G'$ into $\hat G$ from 3.1.3 and the notation for $\hat \g$ from 1.1.2, we can identify $\g'$ as a subalgebra 
of $\hat \g$ as follows
$$
Z = -i(F - G) = -i( \sum_{j=1}^p (U^+_{jj} - U^-_{jj}) - \sum_{j=p+1}^{p+q} (U^+_{jj} - U^-_{jj}) )
$$
and its realization by the infinitesimal Weil \re $d\omega $
$$
\hat Z = \pi (x,x) - (1/4\pi) \Delta ,
$$
where we use the notation indicating the quadratic form given by $S_0$ and its Laplacian
$$
(x,x) := (x,x)_{S_0} :=  \sum_{j=1}^p x_j^2 - \sum_{j=p+1}^{p+q} x_j^2,\,\, \Delta := \Delta _{S_0} := \sum_{j=1}^p \partial _j^2 - \sum_{j=p+1}^{p+q} \partial _j^2.
$$ 
The same way, we have 
$$
X_{\pm} := (1/2)(H \pm i(F+G)) = (1/2)(\sum_{j=1}^n A_{jj} \pm i( \sum_{j=1}^p (U^+_{jj} + U^-_{jj}) - \sum_{j=p+1}^{p+q} (U^+_{jj} + U^-_{jj}) ))
$$
and its realization as a differential operator acting on the Schwartz space $\Ss(\R^n)$
$$
\hat X_{\pm} = (1/2)( E + n/2 \mp (\pi(x,x) + (1/4\pi)\Delta)),
$$
where we use the Euler operator
$$
 E := \sum_{j=1}^n x_j\partial _j.
$$
{\bf 3.3.3} The orthogonal group $G = \Or(S_0) = \Or(p,q)$ has as its Lie algebra
\begin{eqnarray*}
\ot(p,q) = \{ Y =  \begin{pmatrix} Y^{11}&Y^{12} \\Y^{21}&Y^{22} \end{pmatrix}; \,\, 
Y^{11} & = & - \ta Y^{11} \in \M_p(\R),\\ \,Y^{22} & = & - \ta Y^{22} \in \M_q(\R),\,Y^{12} = \ta Y^{21} \in \M_{pq}(\R) \}.
\end{eqnarray*}
One has $\dim \ot(p,q) = n(n-1)/2.$ We write 
$$ 
\ot(p,q) = \ka + \p;\,\, \ka = \{\begin{pmatrix} Y^{11}&\\&Y^{22} \end{pmatrix}\} \simeq \ot(p) \times \ot(q), \,\,
 \p = \{\begin{pmatrix} &Y^{12} \\\ta Y^{12}& \end{pmatrix}\}.
$$ 
As in 3.1.4 $\alpha ,\beta $ denote indices $1,\dots, p\,\,$ and $\mu ,\nu $ indices between $p+1$ and $p+q.$ 
Then $\ka$ is spanned by $n\times n-$matrices of the types
$$
E_{\alpha \beta } - E_{\beta \alpha },\,\, E_{\mu \nu } - E_{\nu \mu }
$$
and $\p$ by those of the type
$$
E_{\alpha \mu } + E_{\mu \alpha }.
$$ 
The embedding of $G$ into $\hat G$ from 3.1.3 induces an embedding of $\g$ into $\hat \g$ given by
$$
\g \ni Y \longmapsto \begin{pmatrix} -\ta Y& \\&Y \end{pmatrix} \in \hat \g.
$$
We use this for an identification and hence with the notation from 1.1.2 can realize the elements of $\ka$ 
in the Weil \re as operators acting on $\Ss(\R^n)$ by
$$
\hat A_{\alpha \beta } - \hat A_{\beta \alpha } = x_\alpha \partial _\beta - x_\beta \partial _\alpha ,\,\,
\hat A_{\mu \nu } - \hat A_{\nu \mu } = x_\mu \partial _\nu - x_\nu \partial _\mu 
$$
and the elements of $\p$ by
$$
\hat A_{\alpha \mu  } + \hat A_{\mu  \alpha } = x_\alpha \partial _\mu  + x_\mu  \partial _\alpha .
$$

{\bf 3.3.4} \Exa $G = \Or(2,1):$ To simplify things, we look at the example $p=2, q=1,$ i.e.~$\g = \ot(2,1) \simeq  \s(2,\R)$. 
Here we use the notation
$$
H := \begin{pmatrix} &1& \\-1&&\\&&0 \end{pmatrix},\,\,Y_1 := \begin{pmatrix} &&1 \\&&&\\1&& \end{pmatrix} ,\,\,
Y_2 =  \begin{pmatrix} &&& \\&&1\\&1& \end{pmatrix}
$$
and have
$$
[H,Y_1] = - Y_2, \,\,[H,Y_2] = Y_1, \,\,[Y_1,Y_2] = H.
$$
With the identification given by the embedding from 3.1.2 above one has
$$
H = A_{12} - A_{21},\,\, Y_1 = - (A_{13} + A_{31}),\,\,Y_2 = - (A_{23} + A_{32})
$$
and the realization as operators for the Weil \re
$$
\hat H = x_1\partial _2 - x_2\partial _1,\,\,\hat Y_1 = - (x_1\partial _3 + x_3\partial _1),\,\,\hat Y_2 = - (x_2\partial _3 + x_3\partial _2).
$$
As usual, we complexify
$$
\g_c = <H_0 := -2iH, \,\,Y_{\pm} := Y_1 \pm iY_2>
$$
with
$$
[H_0,Y_{\pm}] = \pm 2Y_{\pm},\,\,[Y_+,Y_-] = H_0
$$
and get the corresponding operators
$$
\hat H_0 = - 2i(x_1\partial _2 - x_2\partial _1),\,\,\hat Y_{\pm} = - (x_1 \pm ix_2)\partial _3 - x_3(\partial _1 \pm i\partial _2).
$$
{\bf 3.3.5} \Rem: If we apply these operators to the vacuum vector of the Weil \re $\omega ,$ the Gaussian 
$$
\varphi _0(x) = e^{- \pi  (x_1^2 + x_2^2 + x_3^2)},
$$
we get
$$
\hat H_0 \varphi _0 = 0,\,\,\hat Y_{\pm}\varphi _0 = 4\pi x_3(x_1 \pm ix_2)\varphi _0
$$
and using the operators obtained in 3.3.2 specialized to $S_0 = \begin{pmatrix} 1&& \\&1&\\&&-1 \end{pmatrix}$
$$
\hat Z \varphi _0 = (1/2)\varphi _0,\,\,\hat X_+\varphi _0 = (1 - 2\pi (x_1^2+x_2^2))\varphi _0,\,\,
\hat X_-\varphi _0 = ((1/2) - 2\pi x_3^2)\varphi _0.
$$
Hence, for the restricted \re $\omega \vert_{G \times G'},$ the Schwartz function $\varphi _0$ generating the Siegel theta 
is a vector of weight $(0,1/2)$ but not a lowest weight vector as it is for the Weil \re of the ambient group $\hat G.$ 
It is a natural task to search for a Schwartz function which is a vector of dominant weight for irreducible \res 
contained as sub\res in $\omega \vert_{G \times G'}.$ Before we go into this, we just state the following observation as a 
byproduct of the small calculations leading to the Remark above.\\

\subsection {Intermezzo: The Gaussian and $U(\g_c)-$Modules} 

We stay with the example $G = \Or(2,1)$ though a generalization should be easy.\\

{\bf 3.4.1} \Rem: If we apply the operators for the derived \re of the restriction $\omega \vert_{G \times G'}$ to the vector $\varphi _1$ with
$$
\varphi _1(x) := e^{- \pi  (x_1^2 + x_2^2 - x_3^2)} = e^{ - \pi  (x,x)}
$$
we get
$$
\hat H_0 \varphi _1 = 0,\,\,\hat Y_{\pm}\varphi _1 = 0
$$
and
$$
\hat Z \varphi _1 = (3/2)\varphi _1 ,\,\,\hat X_+\varphi _1 = (1/2)(3 - 4\pi (x,x))\varphi _1,\,\,\hat X_-\varphi _1 = 0.
$$
Hence, $\varphi _1,$ which obviously is not a Schwartz function, has the properties of a lowest weight vector of weight 
$0$ for $\omega \vert_G$ and weight $3/2$ for $\omega \vert_{G'}.$ \\
One has
$$
\varphi _1(x) = e^{2\pi x_3^2}\varphi _0(x)
$$
and (from 1.1.3) $\hat U^+_{33} = \pi i x_j^2.$ Then the Remark above also reflects in the formal calculation where
we identify the elements of $\g'_c$ with their images given by the embedding into $\hat \g$: \\

{\bf 3.4.2} \Prop: Let $v_0$ be an element of an $U(\hat \g_c)-$module such that
$$
H_j v_0 = (1/2)v_0 \fa j = 1,\dots,n \,\,{\rm and} \,\,\check U^- v_0 = 0 \,\fa \check U^- \in \syp^{(0,2)},
$$
and 
$$
S := \sum_{l=1}^\infty (-2iU^+_{33})^l / l!.
$$
Then one has
$$
ZSv_0 = (3/2)Sv_0,\,\, X_-Sv_0 = 0.
$$
\pr: We recall
$$
Z = -i(\sum_j \varepsilon _j(U^+_{jj} - U^-_{jj})),\,\, \varepsilon _1 = \varepsilon _2 = - \varepsilon _3 = 1,
$$
$$
X_{\pm} = -(1/2)( \sum_j A_{jj} \pm i \sum_j \varepsilon _j(U^+_{jj} + U^-_{jj}))
$$
and abbreviate
$$
U^+_{33} =: V, \,\,U^-_{33} =: U, \,\,A_{33} =: A.
$$
$V$ commutes with all $U^+_{jj}$ and with $U^-_{jj}$ for $j=1,2$ and one has $[V,U] = A.$ 
Moreover, $V$ commutes with $A_{jj}, j=1,2$ and one has $[V,A] = -2V.$ By induction, one easily verifies for $l \in \N$ 
as relations in $U(\g_c)$
$$
AV^l = V^lA + 2lV^l \,\,{\rm and}\,\, UV^l = V^lU - lV^{l-1}A - l(l-1)V^{l-1}.
$$
Hence we get
$$
AS = SA + \sum 2l(-2i)^lV^l/l! = SA -4iSV \,\,{\rm and}\,\, US = SU + 2iSA + 4SV.
$$
and
\begin{eqnarray*}
X_-Sv_0 & = &  (1/2) (S(\sum A_{jj} - 4iSV - i(\sum \varepsilon _j(U^+_{jj} + U^-_{jj}) - 2iSA - 4SV)v_0\\
         & = &  (1/2) (S(\sum \varepsilon _j(A_{jj} - i(U^+_{jj} + U^-_{jj})v_0 = 0,
\end{eqnarray*}
as $v_0$ has the property $\check U^-{jj} v_0 = 0$ for $j=1,2,3.$ Similarly, one has
$$
ZS = -i(S(\sum \varepsilon _j(U^+_{jj} - U^-_{jj})) + 2iSA + 4SV).
$$
Here we use that we have the relations $-i(U^+_{jj} - U^-_{jj})v_0 = 1/2$ for all $j$ and  $\check U^-_{jj} v_0 = (iA_{jj} + U^+_{jj} +U^-_{jj}) v_0 = 0$ 
leading to $(iA_{jj} + 2U^+_{jj})v_0 = (i/2)v_0$ and get
$$
ZSv_0 = (3/2)v_0.
$$

{\bf 3.4.3} From the Remark 3.4.1 above one would expect to have also $H_0Sv_0 = Y_\pm v_0 = 0.$ Here 
again we identify $\g = \ot(2,1)$ with its image in $\hat \g_c.$ One has 
$$
H_0 = -2i H, \,\, Y_{\pm} = Y_1 \pm iY_2
$$
with
$$
H = A_{12} - A_{21}, \,\, Y_1 = - (A_{13} + A_{31}, \,\,Y_2 = - (A_{23} + A_{32}).
$$
$H$ commutes with $V = U^+_{33}.$ One has
$$
[Y_1,U^+_{33}] = - U^+_{13}, \,\,[Y_2,U^+_{33}] = - U^+_{23}, 
$$
and hence
$$
Y_1V^l = V^lY_1 - lv^{l-1}U^+_{13},\,\,Y_2V^l = V^lY_2 - lv^{l-1}U^+_{23}.
$$
We get
$$
H_0S = SH_0,\,\, {\rm and} \,\, Y_{\pm}S = S(Y_{\pm} + 2i(U^+_{13} \pm iU^+_{23}))
$$
and see that for the relations $H_0Sv_0 = 0, \, Y_{\pm}Sv_0 = 0$ one needs the conditions
$$
A_{12}v_0 = A_{21}v_0,\,,\,\, {\rm and} \,\, (A_{jk} + A_{kj})v_0 = 2iU^+_{jk}v_0.
$$
These conditions are fulfilled if we take $v_0 = \varphi _0$ and the realization by the Weil \re but certainly there is 
more background to this simple discussion. \\

{\bf 3.4.4} It should be interesting to see an explicit decomposition of the $U(\hat \g_c)$ module belonging to the 
Weil \re into its irreducible $(U(\g_c) \times U(\g'_c))-$modules. I don't know whether this is done someplace. 
But a somewhat equivalent task is easily accessible and we will describe this now.\\

\section  [Weil Representation and Rallis-Schiffmann Thetas]
{Decomposition of the Weil Representation and Rallis-Schiffmann Thetas}

\subsection {Weil Representation $\tilde \omega _S$}

 We present explicit material (partially going back to Rallis and Schiffmann [RS1-3]) collected by M. Vergne in [LV] 
concerning the decomposition of the Weil \re $\omega $ as a \re of the dual pair $G = \Or(p,q), G' = \Spe(m,\R),$ in 
particular for $m=1$. In this section we are interested in the discrete spectrum of $\omega $ as 
\re of $G\cdot G'.$ It is a special case of a more general conjecture by Howe that the restriction of $\omega $ 
induces a one-one correspondence between the irreducible unitary (discrete series) \res of $G$ and $G'.$ This has been proved by 
work of Howe, Rallis, Schiffmann, and Strichartz ([Ho3], [RS1], [Str]). We shall follow essentially the presentation given in [LV] p.205ff.\\

{\bf 4.1.1} We restrict our treatment to the case $G' = \SL(2,\R), G = \Or(S) = \Or(p,q)$ and denote by $\tilde \omega _S$ the restriction to $G\cdot \tilde G'$ of 
the Weil \re $\tilde \omega $ of the metaplectic cover $\tilde {\hat G},  \hat G = \Spe(n,\R), n = p+q.$ In order to describe the 
decomposition of $\tilde \omega $ we have to fix a lot of notation and reproduce some elements of the \re theory of $G'$ and $G$ 
though the reader will probably know most of this.\\

\subsection {Discrete Series of $G' = \SL(2,\R)$} 

This group has three types of unitary irreducible \res, the discrete, the principal, and the complementary series \res. \\

{\bf 4.2.1} Here we only report on the discrete series $\pi_k^+$ and $\pi_k^-$ for $k \in \N, k\geq 2$. $\pi_k^+$ is a lowest weight \re with a 
lowest weight vector $v$ which, using our notation from 3.6, is characterized by  
$$
d\pi_k^+(Z) v = k v,\,\, d\pi_k^+(X_-) v = 0.
$$
Similarly $\pi_k^-$ has a vector $v_-$ of highest weight $-k$ characterized by 
$$
d\pi_k^-(Z) v_- = k v_-, d\pi_k^-(X_+) v_- = 0.
$$
Standard models for these \res come from the space $\OO(\Ha)$ of holomorphic functions on the 
upper half plane, resp.~the space $\bar{ \OO}(\Ha)$ of antiholomorphic functions with prescriptions
$$
(\pi_k(g^{-1})f) = (c\tau + d)^{-k}f(g(\tau )); \fa f \in \OO(\Ha), \,\, g = \begin{pmatrix} a&b\\c&d \end{pmatrix}\,\,\in G',
$$
resp.
$$
(\bar \pi_k(g^{-1})f) = (c\bar \tau + d)^{-k}f(g(\tau )); \fa f \in \bar {\OO}(\Ha), \,\, g = \begin{pmatrix} a&b\\c&d \end{pmatrix}\,\,\in G'.
$$
If one restricts to the spaces $L^2_{hol}(\Ha,d\mu _k)$ resp.~$L^2_{antihol}(\Ha,d\mu _k)$ of holomorphic 
resp.~antiholomorphic funtions on $\Ha$ with finite norm for the measure $d\mu _k(\tau ) = v^{k-2}dudv,$ one has unitary \res. \\
 
{\bf 4.2.2} \Rem: It is an easy but interesting exercise to verify that 
$$
\psi _w(\tau ) := (\bar \tau - w)^{-k}, \,\, w \in \Ha,
$$
is a lowest weight vector for $\bar \pi_k$ which fulfills the {\it fundamental relation}
$$
\bar \pi_k(g) \psi _w = (cw + d)^{-k}\psi _{g(w)}.
$$
We use the coordinization
$$
G' = \SL(2,\R) \ni g = \begin{pmatrix} a&b\\c&d \end{pmatrix}\,= n(u)t(v)r(\vartheta ),
$$
with
$$
n(u) := \begin{pmatrix} 1&u\\&1 \end{pmatrix},\,\,t(v) := \begin{pmatrix} v^{1/2}&\\&v^{-1/2} \end{pmatrix},\,\,
r(\vartheta ) := \begin{pmatrix} \cos \vartheta &\sin \vartheta \\-\sin \vartheta &\cos \vartheta  \end{pmatrix}
$$
and the standard notation $j(g,\tau ) = (c\tau + d),\,\, j_k(g,\tau ) = (c\tau + d)^{-k}.$
Hence, we have $g(i) = u +iv = \tau, j(g,i) = (ci+d) = v^{-1/2}e^{-i\vartheta } $ and we put 
$$
g_\tau := n(u)t(v).
$$
(This coordinization goes back to Lang's book $\SL_2(\R).$ Vergne in [LV] and some other authors use $u(\theta ) := r(-\theta )$
 for the parametrization of $\SO(2)$)\\

For $U \in \g'$ we denote by $\RR_U$ the right invariant operator on $\Cl^\infty (G')$ given by
$$
\RR_U \Phi (g) = \frac{d}{dt} \Phi (exp(-tU)g)\vert_{t=0}
$$
and by $\LL_U$ the left invariant operator given by
$$
\LL_U \Phi (g) = \frac{d}{dt} \Phi (g(exp(tU)))\vert_{t=0}.
$$
{\bf 4.2.3} For a moment we replace $u,v$ by $x,y.$ A standard calculation (see for instance Lang [La] p.113) leads to
$$
\begin{array}{rcl}
\RR_Z & = & i((1+x^2-y^2) \partial _x + 2xy\partial _y + y\partial _\vartheta) \\
\RR_{X_{\pm}} & = & \pm (i/2) (((x \pm i)^2 - y^2)\partial _x + 2(x \pm i)y\partial _y + y\partial _\vartheta).
\end{array}
$$
and
$$
\begin{array}{rcl}
\label{lld}
\LL_Z & = & - i \partial _\vartheta \\
\LL_{X_{\pm}} & = & \pm (i/2) e^{\pm 2i\vartheta }(2y(\partial _x \mp i \partial _y) - \partial _\vartheta).
\end{array}
$$

If $v_0 \in V$ is a lowest weight vector of weight $k$ for the \re $(\pi,V),$ one has $\pi(r(\vartheta ))v_0 = e^{ik\vartheta }v_0$ 
and one can see that
$$
(ci+d)^k\pi(g)v_0 = v^{-k/2}\pi(g_\tau )v_0
$$
depends only on $\tau ,$ hence we abbreviate it by $v_\tau .$ \\

{\bf 4.2.4} \Rem: The lowest weight \re $(\pi,V)$ is isomorphic to a sub\re of $(\bar \pi_k, \bar {\OO}(\Ha)),$
the isomorphism being obtained by sending $v_w$ to $\psi _w$.\\

If one has a group acting transitively on a space and an automorphic factor $j$ for this action, there is a standard procedure 
to define an appropriate lift of functions on the space to functions living on the group. 
In our case, we look at the action of $G'$ on $\Ha$ and the automorphic factor $j_k(g,\tau ) = (c\tau + d)^{-k}$ and for a 
function $f : \Ha \longrightarrow \C$ define the lift $\varphi _k: f \longmapsto \Phi _f $ by $\Phi _f(g) := j_k(g,i) f(g(i)).$ 
One has the fundamental fact (see for instance [La] IX \S 5):\\

{\bf 4.2.5} \Prop: The function $f$ is holomorphic if and only if
$$
\LL_{X_-} \Phi _f = 0.
$$

We denote
$$
M(G',k) := \{ \Phi \in \Cl(G'); \Phi (gr(\vartheta )) = e^{ik\vartheta }\Phi (g) \fa g \in G',\,r(\vartheta ) \in K'\}
$$
and by $\lambda _k$ the \re of $G'$ given on $M(G',k)$ by (inverse) left translation $\lambda _k(g_0)\Phi (g) = \Phi (g_0^{-1}g).$ 
Then it is not difficult to prove \\

{\bf 4.2.6} \Prop: The lifting map $\varphi _k$ intertwines the \res $\pi_k$ and $\lambda _k.$\\

There is the inverse map $I_k$ to $\varphi _k$ given by associating to $\Phi \in M(G',k)$ the function $f(\tau ) = v^{-k/2}\Phi (g_\tau ).$\\

{\bf 4.2.7} Up to now the index $k$ was an integer and our group $G' = \SL(2,\R)$. We also will have to work with the metaplectic 
cover $\tilde G' = \Mp(2,\R)$ and halfinteger $k$. As in this report we follow Vergne in [LV] p.184, we look at the universal cover $\GG$ of 
$G'$ in the form
$$
\GG := \{ (g, \varphi _g); g \in G', \varphi \in \OO(\Ha) \,\,{\rm such\,\,that}\,\, e^{\varphi _g(\tau )} = j(g,\tau ) \}.
$$
We denote by $pr$ the projection $\GG \longrightarrow G'.$ and, for $\alpha \in \R,$ by $\pi_{\alpha ,0}$ the \re of $\GG$ given on $\OO(\Ha)$ 
by
$$
\pi_{\alpha ,0}((g,\varphi )^{-1}f(\tau ) = e^{-\alpha \varphi (\tau )}f(g(\tau )).
$$
For $\alpha = k \in \Z$ the projection $pr$ intertwines $\pi_{k,0}$ with the \re $\pi_k$ of $G'$ given above. \\

For $J := \begin{pmatrix} &1\\-1& \end{pmatrix}$ one has the one parameter subgroup 
$$
K' = \SO(2) = \{ r(\vartheta ) = \begin{pmatrix} \cos \vartheta &\sin \vartheta \\-\sin \vartheta &\cos \vartheta  \end{pmatrix} \}
$$
of $G'$ and, above in $\GG,$ the group 
$$
\KK = \{ \delta (\vartheta ) := (r(\vartheta ), \varphi _\vartheta ); \varphi _\vartheta (i) = -i\vartheta  \}.
$$ 
For $\alpha \in \R$ we define
$$
M(\GG,\alpha ) := \{ \Phi \,\,{\rm analytic \,\,on}\,\, \GG; \,\,\Phi ({\tilde g} \delta (\vartheta )) = e^{i \alpha \vartheta } \Phi ({\tilde g}) \fa {\tilde g} \in \GG  \}.
$$
and denote by $\lambda _\alpha $ the \re of $\GG$ given on $M(\GG,\alpha )$ by  left inverse multiplication 
$\lambda _\alpha (\tilde g_0)\Phi  = \Phi (\tilde g_0^{-1}\tilde g).$
We put $a_\alpha ((g,\varphi ) := e^{\alpha \varphi (i)}.$ One can verify easily that for $\Phi \in M(\GG,\alpha )$ the function $I_\alpha \Phi $ 
given by
$$
I_\alpha (\tilde g) := a_\alpha (\tilde g)\Phi (\tilde g)
$$
is invariant by right translation by $\delta (\vartheta ),$ Hence one can find a function on $\Ha$ still denoted by $I_\alpha \Phi $ 
such that $I_\alpha \Phi (\tilde g) = I_\alpha \Phi (pr(\tilde g))(i)).$ From Lemma 2.3.13 in [LV] we take over the following fact.\\

{\bf 4.2.8} \Rem: $I_\alpha $ intertwines the \res $\lambda _\alpha $ and $\pi_{\alpha ,0}$\\

\subsection {Elementary Thetas and their Construction}

 Though it is somewhat redundant in view of what we already did in 1.2.2, we want to illustrate the notions we just introduced and 
we describe again a bit more precisely a construction of the simplest theta functions. \\

{\bf 4.3.1} The Weil \re of $G'$ on $\Hi = L^2(\R)$ as a special case given by the formulae in 1.1.1 decomposes into two irreducible 
sub\res $\omega _{even} =: \pi_{W,+}, \omega _{odd} =: \pi_{W,-}$ namely those given by the even and the odd functions. 
As already shown in 1.1.4, the even function $\varphi _0 = e^{- \pi x^2}$ is a lowest weight vector of  weight $1/2$ for $\omega _{even}$ and 
a small calculation shows that 
$$
\varphi '_0(x) := xe^{- \pi x^2}
$$
is a lowest weight vector of weight $3/2$ for $\omega _{odd}.$ \\

{\bf 4.3.2 Vergne's Program:} As Vergne states in [LV] p.179, one can construct a modular form of weight $k$ and character $\chi $ 
for a discrete subgroup $\Gamma $ of $G'$ (or $\GG$) 
by achieving two steps.\\
\phantom{}\hspace{.5cm}1.) Find lowest a lowest weight vector $v \not= 0$ of weight $k$ in the space $\Hi$ of a \re $\pi.$\\
\phantom{}\hspace{.5cm}2.) Find a functional $\theta \in \Hi'$ such that 
$$
\pi(\gamma )\theta = (\chi (\gamma ))^{-1}\theta \fa \gamma \in \Gamma .
$$
Then, using the map $I_k$ resp $I_\alpha $ introduced above, one can transfer $\theta (\pi(g)v)$ to 
a holomorphic function $f$ living on $\Ha$ with the properties wanted.\\

In our case we have the lowest weight vectors $v = \varphi _0$ and $v = \varphi '_0$  of weights $k = 1/2$ and $k = 3/2$ 
for the sub\res of the Weil \re. For the determination of the semi invariant functional fitting to certain subgroups Vergne 
goes a long way back to the Schr\"odinger and lattice \re of the Heisenberg group. To abbreviate things, 
here we simply take the theta distribution $\theta = \sum_{n\in\Z} \delta _n$ and get
$$
\theta (\omega( g_\tau )\varphi _0) = v^{1/4} \sum_{n\in\Z} e^{\pi i\tau n^2},
\, \, I_{1/2}\theta (\omega (g_\tau )\varphi _0) = \sum_{n\in\Z} e^{\pi i\tau n^2} = \,\theta (\tau ) 
$$
and
$$
\theta (\omega (g_\tau )\varphi' _0) = v^{3/4} \sum_{n\in\Z} ne^{\pi i\tau n^2},
\, \, I_{3/2}\theta (\omega (g_\tau )\varphi' _0) = \sum_{n\in\Z} e^{\pi i\tau n^2} .
$$
The first function comes out as our fundamental theta function with the precise transformation property given as follows 
(see for instance [LV] p.204)\\
{\bf 4.3.3} \Th: For every $\gamma = \begin{pmatrix} a&b\\c&d \end{pmatrix}$ with $a,b,c,d \in \Z,\, ad-bc =1, \,ac\equiv 0 \mod 2,\, bd\equiv 0 \mod 2$ 
one has
$$
\theta (\gamma (\tau )) = \lambda (\gamma )(c\tau + d)^{1/2}\theta (\tau )
$$
where $(c\tau + d)^{1/2}$ is the principal determination of $(c\tau + d)^{1/2}$ and $\lambda $ is a rather complicated character 
given in [LV] p.201. In particular, if $c$ is even, $\lambda (\gamma ) = \varepsilon _d^{-1}(\frac{2c}{d}).$ 
For $d$ an odd integer $\varepsilon _d$ is defined to be
$$
\begin{array}{rcl}
\varepsilon _d & = & 1 \,\,{\rm if}\,\, d \equiv 1 \mod 4\\
& = & i \,\,{\rm if}\,\, d \equiv -1 \mod 4.
\end{array}
$$
The second function apparently is identically zero. But one still gets reasonable functions by taking a character $\psi \mod N$ and $t\in\Z$
defining 
$$
\theta ^+_{\psi,t} (\tau ) := \sum_{n\in\Z} \psi (n) e^{2\pi i\tau tn^2}\,\,{\rm for}\,\, \psi \,\,{\rm even}
$$
and
$$
\theta ^-_\psi {(\tau,t} ) := \sum_{n\in\Z} \psi (n) ne^{2\pi i\tau tn^2}\,\,{\rm for}\,\, \psi \,\,{\rm odd}
$$
with
$$
\theta ^+_{\psi,t} (\gamma (\tau )) = \psi (d)(\frac{t}{d}) \varepsilon ^{-1}_d(c\tau + d)^{1/2}\theta ^+_{\psi,t} (\tau )\,\fa \gamma \in \Gamma _0(4N^2t),
$$
and
$$
\theta ^-_{\psi,t} (\gamma (\tau )) = \psi (d)(\frac{t}{d}) \varepsilon ^{-1}_d(c\tau + d)^{3/2}\theta ^-_{\psi,t} (\tau )\,\fa \gamma \in \Gamma _0(4N^2t).
$$
\subsection {The Theta Miracle} 

We add an observation which for me plays the role of a {\it Theta Miracle}. In 1.3.2 and 1.3.3 we extended the Weil \re of $G'$ resp.~$\GG$ 
with the Schr\"odinger \re of the Heisenberg group to the 
\re $\pi_{SW}$ of the Jacobi group $G^J$. If we take our vacuum vector of lowest weight $1/2,$ namely $\varphi _0(x) = e^{-\pi x^2},$ 
we get 
$$\varphi _{g^J}(x) := 
\pi_{SW}((p,q,k) g_\tau r(\vartheta ))\varphi _0(x) = v^{1/4} e^{\pi i\vartheta /2}e^{\pi i(\kappa + p^2\tau  + pq + x^2\tau + 2(p\tau + q)x)}.
$$
We will see that $\Phi _0(g^J) := \varphi _{g^J}(0)$ directly can be interpreted as a lowest weight vector. \\

{\bf 4.4.1} To explain this and prepare some more material useful later, we treat the Lie algebra of the Jacobi group.
As in [BeS] or [Ya], we describe the Lie algebra
${\frak g}^J$ as a subalgebra of ${\frak s}{\frak p}(2,\bf{R})$ by
$$
G(x,y,z,p,q,r)=
\begin{pmatrix}
x&0&y&q\\p&0&q&r\\z&0&-x&-p\\0&0&0&0
\end{pmatrix}
$$
and denote

$$
X=G(1,0,\dots,0),\dots,R=G(0,\dots,0,1).
$$

We get the commutators

$$
\begin{array}{lll}
[X,Y]= \, 2Y,\,&[X,Z]=-2Z,\,&[Y,\,Z]=X,\\ {[}X,\,P]= - P,\,&[X,Q]= \,\,Q,\,&[P,Q]=2R,\\{[}\,Y,P]=-Q,\,&[\,Z,Q]= - P,&
\end{array}
$$

all others are zero. Hence, we have the complexified Lie algebra given by 
$$
{\frak g}_c^J = <Z_1,X_\pm,Y_\pm,Z_0>
$$
where as in [BeS]
$$
\begin{array}{l}
Z_1=-i(Y-Z),\ Z_0=-iR,\\
\\
X_\pm=(1/2)(X \,\pm\, i(Y+Z)),\ Y_\pm=(1/2)(P\, \pm\, iQ)
\end{array}
$$
with the commutation relations
$$
[Z_1 , X_{\pm}] = {\pm}2X_{\pm} ,\,  [Z_0 , Y_{\pm}] = {\pm}Y_{\pm} ,  \,  {\rm etc}.
 $$
{\bf 4.4.2} The complexified
Lie algebra ${\frak g}_c^J$ of the Jacobi group
is realized (see [BeS] p.12) by the left invariant differential operators
$$
\begin{array}{l}
{\mathcal L}_{Z_0}=i\partial_\kappa\\
\medskip
{\mathcal L}_{Y_\pm}=(1/2)y^{-1/2}e^{\pm i\vartheta}
(\partial_p-(x\pm iy)\partial_q-(p(x+iy)+q)\partial_\kappa)\\
\medskip
{\mathcal L}_{X_\pm}=\pm(i/2)e^{\pm2i\vartheta}(2y
(\partial_x\mp i\partial_y)-\partial_\vartheta)\\
\medskip
{\mathcal L}_Z=-i\partial_\vartheta
\end{array}
$$
acting on differentiable functions $\phi=\phi(g)$ with the coordinates 
coming from
$$
g=(p,q,\kappa)n(x)t(y)r(\vartheta).
$$
where
$$
n(x)=\left(\begin{array}{cc}1&x\\&1\end{array}\right),t(a)=
\left(\begin{array}{cc}a^{1/2}&\\&a^{-1/2}\end{array}\right),\ r(\vartheta)=
\left(\begin{array}{rc}\alpha&\beta\\-\beta&\alpha\end{array}\right)
$$
with
$$
 \alpha=\cos\vartheta,\ \beta=\sin\vartheta,
$$
and 
$$
g(i,0)=(\tau,z)=(x+iy,p\tau+q).
$$
As usual, we put
$$
N' = \{n(x), x \in {\bf R}\}, A' = \{t(a), a \in {\bf R}_{>0}\},
$$ 
$$
K' = {\rm SO}(2) = \{r(\vartheta), \vartheta \in {\bf R}\}, M = \{{\pm} E\}.
$$
{\bf 4.4.3} We introduce the standard automorphic factor for $g^J = n(u)t(v)r(\vartheta )(\lambda ,\mu ,\kappa ) \in G^J$
$$
j_{m,k}(g^J,(\tau ,z)) := (c\tau + d)^{-k}e^m(\kappa - \frac{c(z+\lambda \tau + \mu )^2}{c\tau +d} + \lambda ^2\tau + 2\lambda z + \lambda \mu )
$$
{\bf 4.4.4} \Prop: The function
$$
\Phi _0(g^J) = v^{1/4} e^{i\vartheta /2}e^{\pi i(\kappa + pz)},\,\, g^J = (p,q,\kappa ) n(u)t(v)r(\vartheta )
$$
fulfills the relations
$$
\LL_{Z_0}\Phi _0 = (1/2) \Phi _0, \,\,\LL_{Z}\Phi _0 = (1/2) \Phi _0,\,\, \LL_{Y_-}\Phi _0 = 0,\,\, \LL_{X_-}\Phi _0 = 0,
$$
i.e.~$\Phi _0$ spans in the restriction of the right regular \re of $G^J$ a \re of lowest weight $(1/2,1/2)$ for the subgroup 
$K^J := \{(0,0,\kappa )r(\vartheta ); \kappa ,\vartheta \in \R \} $ of $G^J.$\\

{\bf 4.4.5} It is easy to see that one has $\Phi _0(g^J) = j_{1/2,1/2}(g^J,(i,0))$ and that $\Phi _0$ is invariant under left multiplication 
by elements of the group $N^J := \{(0,q,0 )n(u); q,u \in \R\}.$  
Now, if we want invariance under the group $\Gamma ^J := \SL(2,\Z)\ltimes Z^3$ or an appropriate subgroup, 
we try the averaging
$$
\Theta (g^J) := \sum_{\ell \in \Z}\Phi ((\ell,0,0)g^J),
$$
which has the same outcome as the application of the theta distribution from 4.3.2 to 
$\varphi _{g^J}(x):$ We get
$$
\begin{array}{rcl}
\Theta (g^J) & = & \sum_{\ell\in\Z} \pi_{SW}((p,q,k) g_\tau r(\vartheta ))\varphi _0(\ell) \\*[0.2cm]
& = & v^{1/4} e^{\pi i\vartheta /2}e^{\pi i(\kappa + p^2\tau  + pq)}\sum_{\ell\in\Z} e^{\pi i(\tau \ell^2+ 2(p\tau + q)\ell)}\\*[0.2cm]

& = & j_{1/2,1/2}(g^J,i)\theta (\tau ,z).
\end{array}
$$
As usual, using the Poisson summation formula, one can see ([Be1]) that (by the {\it theta miracle}\,) 
our construction has a built in invariance property 
with respect to the element $(0,0,0)J \in \Gamma ^J$ and hence for the whole theta group. We cite from [MuI] p.32.\\

{\bf 4.4.6} \Th: For $\gamma = \begin{pmatrix}a&b\\c&d\end{pmatrix} \in\Gamma _\theta $ one has
$$
\theta (\frac{a\tau + b}{c\tau + d},\frac{z}{c\tau + d}) = \zeta (c\tau + d)^{1/2}e^{\pi i cz^2/(c\tau + d)}\theta (\tau ,z)
$$
with
$$
\begin{array}{rcll}
\zeta & = & i^{(d-1)/2}(\frac{c}{\vert d\vert}) \,\,&c\,\,{\rm even},\,\,d\,\,{\rm odd}\\
& = & e^{\pi ic/4}(\frac{d}{c}) \,\,&d\,\,{\rm even},\,\,c\,\,{\rm odd}.
\end{array}
$$
{\bf 4.4.7} From [MuI] p.227 we take over that one has the same relation for the more general function
$$
\theta \, [\begin{array}{c}a\\b\end{array}](\tau ,z) := \sum_{n\in\Z} e^{\pi i((n + a)^2\tau + 2(n +a)(z + b))}.
$$
If one takes here the differentiation $\partial _z$ and afterwards puts $z = 0,$ one gets
$$
\frac{\partial \theta\, [\begin{array}{c}a\\b\end{array}]}{\partial z}(\gamma( \tau) ,0) = 
\zeta (c\tau + d)^{3/2}\frac{\partial \theta [\begin{array}{c}a\\b\end{array}]}{\partial z}( \tau ,0)
$$
i.e.~
$$
\frac{\partial \theta\, [\begin{array}{c}a\\b\end{array}]}{\partial z}(\tau ,0) = 2\pi i\sum_{n\in\Z}(n + a)e^{\pi i[(n + a)^2\tau + 2(n + a)b)}
$$
is a modular form of weight $3/2.$ We observe that the function $\varphi '_0$ reappears by a differentiation of a function 
stemming from $\varphi _0.$\\

\subsection {Representations of the Orthogonal Groups}

We start by looking at the compact group $G = \Or(p)$ belonging to a positive definite form $S$ on $V = \R^p.$ 
$G$ acts on functions $f$ on $V$ via $\lambda (g)f(x) := f(g^{-1}x)$ and similarly on the space $\Pp := \C[x_1,\dots,x_p]$ of polynomials. 
One identifies $V$ with its dual $V^*$ via $S.$
and this identification extends to an identification of $\Pp$ the space $\DD$ of differential operators with constant coefficients on $\R^p$ 
via $P[x_1,\dots,x_p] \longmapsto P(\partial _1,\dots,\partial _p).$ Hence $\Pp$ is provided with an $\Or(P)$-invariant hermitian inner product 
given by
$$
<P,Q> := (P(\partial _x)\bar Q)(0).
$$
{\bf 4.5.1 \Prop:} The Laplace operator $\Delta _S := \sum_{j=1}^p \partial _j^2$ commutes with the action of $G$ and the space 
$$
\Hi^m = \Hi^m(\R^p) := \{P \in \Pp; \Delta _S P = 0, \deg P = m \}
$$
of harmonic polynomials of degree $m$ is the space of an irreducible unitary \re $\delta _m$ of $\Or(p).$ \\

One has 
$$
\dim \Hi^m(\R^p) = {m+p-1 \choose p-1} - {m+p-3 \choose p-1}.
$$
If $p=1,$ this is zero for $m\geq 2.$ \\

Moreover, it is a standard fact that each polyomial $P$ can be written as a sum 
$P = \sum S^jP_j$ where the $P_j$ are harmonic and $S = S(x) := S(x,x).$ As the space 
$$
\{ P(x) e^{-\pi S(x)}; \,P \in \Pp \} 
$$ 
is dense in $L^2(\R^p)$, one has a decomposition
$$
L^2(\R^p) = \oplus _{m\in \N_0} L^2(\delta _m) \cong  \oplus _{m\in \N_0} \Hi^m(\R^p)
$$
 where $L^2(\delta _m)$ indicates the isotypic component of type $\delta _m$ in $L^2(\R^p)$ (and $m$ is only $0$ or $1$ for $p=1$). \\

{\bf 4.5.2} If $S^{p-1}$ denotes the sphere in $\R^p,$ one has as well (see for instance [Kn] p.81)
$$
\Cl^\infty(S^{p-1}) \cong \oplus _{m\in \N_0} \Hi^m(\R^p).
$$ 
 Now, we will take a look at the non-compact case $G = \Or(p,q), p\geq q, q\geq 1.$\\
 
{\bf 4.5.3} It is clear that each \re of $G = \Or(p,q)$ will decompose under the compact subgroup $K = \Or(p)\times \Or(q)$ 
and we shall have to find a way to analyse this decomposition.
We take again the action of $g \in G$ given on functions $f$ on $\R^{p+q} = \R^n$ by $\lambda (g)f(w) := f(g^{-1}w).$ We write
$$
\R^n = \R^{p+q} \ni w = (x,y),\,\, x \in \R^p, y \in \R^q;\,\, r(x) := (x^2_1,\dots,x_p^2)^{1/2}, \, r(y) := (y_1^2,\dots,y_q^2)^{1/2}
$$
and
$$
\begin{array}{rcl}
\X^0 & := & \{ (x,y) \in \R^{p+q}; \,\, r(x)^2 - r(y)^2 = 0 \},\\
\X^{\pm} & := & \{ (x,y) \in \R^{p+q}; \,\, r(x)^2 - r(y)^2 \gtrless 0 \},\\
\X^t & := & \{ (x,y) \in \R^{p+q}; \,\, r(x)^2 - r(y)^2 = t^2\}\,\,{\rm for} \,\,t>0.
\end{array}
$$
$G$ acts transitively and one can construct \res on spaces of functions living on each one of these sets. 
Though it is not the most appropriate one for our later application, we briefly follow the nice and accessible 
presentation given in [HT], which is based on the cone $\X^0.$ There are fundamental papers by Strichartz [Str] and 
Rallis-Schiffmann [RS1-3], which use $\X^+$ and the hyperboloid $\X^t.$ We will come to this later but for now 
study spaces of homogeneous functions on the ($G-$invariant) light cone $\X^0:$ For $a \in \C$ denote by $S^a(\X^0)$ the space 
$$
S^a(\X^0) := \{ f \in \C^\infty (\X^0); \, f(tw) = t^af(w), \,w \in \X^0,\, t \in \R_{>0} \}.
$$ 
Since $\Or(p,q)$ commutes with scalar dilatations, it is clear that $S^a(\X^0)$ will be invariant under the action given by $\lambda .$ 
To study the structure of $S^a(\X^0)$ as $\Or(p,q)-$module we consider the action of the compact subgroup $K = \Or(p)\times \Or(q).$ 
We take the tensor product of the harmonic $m-$forms $\Hi^m(\R^p)$ and the $n-$forms $\Hi^n(\R^q)$ and its embedding 
into $S^a(\X^0)$ given by
$$
j_a = j_{a,m,n} :  \Hi^m(\R^p) \otimes \Hi^n(\R^q) \longrightarrow S^a(\X^0); \,\, h_1 \otimes h_2 \longmapsto h_1(x)h_2(y)r(x)^{2b}
$$
with $m+n+2b = a.$ Then one has (Lemma 2.2 in [HT])\\

{\bf 4.5.4 \Le}: As an $\Or(p) \times \Or(q)$ module, the space $S^a(\X^0)$ decomposes into a direct sum
$$
S^a(\X^0) \simeq \sum_{m,n\geq 0} j_a(\Hi^m(\R^p) \otimes \Hi^n(\R^q))
$$
of mutually inequivalent irreducible \res.\\

One refers to the spaces $j_a(\Hi^m(\R^p) \otimes \Hi^n(\R^q))$ as the {\it $K-$types} of $S^a(\X^0).$ To understand the \re 
one has to study how $\Or(p,q)$ transforms one $K-$type into another. For this purpose one can use the {\it ladder operators} 
given by the generators of the Lie algebra $\p$ from 3.3.3
$$
\hat A_{\alpha \mu  } + \hat A_{\mu  \alpha } = x_\alpha \partial _\mu  + x_\mu  \partial _\alpha .
$$
in their action on functions on $\R^{p+q}.$ By a routine calculation (carried out in [HT] \S6.1) this action can be transfered 
to maps between the $K-$types and, for $q>1,$ one gets (Lemma 2.3 in [HT])\\

{\bf 4.5.5 \Le}: For each pair $(m,n)$, there are maps
$$
T^{\pm,\pm}_{m,n}: \p \otimes (\Hi^m(\R^p) \otimes \Hi^n(\R^q)) \longrightarrow (\Hi^{m\pm1}(\R^p) \otimes \Hi^{n\pm1}(\R^q)),
$$
which are independent of $a$ and which are nonzero as long as the target space is nonzero, such that the action of 
$Y \in \p$ on the $K-$type $j_a(\Hi^m(\R^p) \otimes \Hi^n(\R^q))$ is described by the formula
$$
\begin{array}{rcl}
\rho (Y)j_a(\phi ) & = & (a-m-n)j_a(T^{++}_{m,n}(Y \otimes \phi ))\\
&  & + (a-m+n+q-2)j_a(T^{+-}_{m,n}(Y \otimes \phi ))\\
&  & + (a+m-n+p-2)j_a(T^{-+}_{m,n}(Y \otimes \phi ))\\
&  & + (a+m+n+p+q-4)j_a(T^{--}_{m,n}(Y \otimes \phi ),\\
\end{array}
$$
where $\phi \in \Hi^m(\R^p) \otimes \Hi^n(\R^q).$ \\

One sees easily that the {\it transition coefficients} relating the different $K-$types are never zero if $a$ is not an integer 
and hence one has an irreducible $\Or(p,q)$ module. In [HT] one finds a detailed study what happens for integer $a$ and 
which \res are unitary. \\

{\bf 4.5.6} As we are particularly interested in this special case, we finish by reproducing some results concerning the case 
$q = 1$ omitted above. When $q = 1,$ the light cone $\X^0$ is not connected. One denotes
$$
\X^{0\pm} := \{ (x,y) \in \X^0; \,\, y = \pm r(x) \}.
$$
Each subspace $\X^{0\pm}$ is stabilized by a subgroup $\Or^+(p,1)$ of index 2 in $\Or(p,1).$
And analyzing $S^a(\X^0)$ as an $\Or(p,1)$ module is the same as analyzing the space of even functions $S^{a+}(\X^0) \simeq S^a(\X^{0+})$ 
as an $\Or^+(p,1)$ module. For $a\in \C$ we can define embeddings
$$
j_a = j_{a,m} : \Hi^m(\R^p) \longrightarrow S^a(\X^{0+}); \,\, h \longmapsto h(x,y)y^{a-m}.
$$ 
One has 
$$
S^a(\X^{0+}) \simeq \sum_{m \geq 0} j_a(\Hi^m(\R^p))
$$
and parallel to Lemma 4.5.5 the action of $\p$ in this case leads to
$$
(x_j\partial _y + y\partial _{x_j})(hy^{a-m}) = (a-m)T^+_j(h)y^{a-m-1} + (a+m+p-2)T^-_j(h) y^{a-m+1}.
$$
One concludes that $S^a(\X^{0+})$ is always irreducible except when $a$ is an integer, with either $a \geq 0$ or $a \leq -p + 1$, in which case 
$S^a(\X^{0+})$ has two constituents, one of which is finite dimensional. In [HT] there is a discussion of the unitarity of the \res. 
In this degenerate case we have unitarity for $a \in \C$ with Re $a = -(p-1)/2, a \in (-(p-1),0),$ or $a \in \N.$\\

\subsection [Weil Representation Associated to a Definite Quadratic Form]
{Decomposition of the Weil Representation Associated to a Positive Definite Quadratic Form}

We use the notation $(V,S)$ to denote an orthogonal space and $(V',B)$ resp.~$(\hat V := V \otimes V', \hat B := S \otimes B)$ 
to denote symplectic spaces with standard bases as in 3.1.4, $\dim V = n, \dim V' = 2m$. $\omega $ is the Weil \re of 
$\hat G = \Spe(mn,\R) \simeq \Spe(S \otimes  B).$ (resp.~$\tilde \omega $ for the metaplectic cover $\tilde G'$) as in 1.1 and 
$\omega _S$ its restriction to $\Or(S) \times \Spe(B)$ resp.~$\tilde \omega _S$ for the metaplectic cover.
For positive definite $S$ the decomposition of $\tilde \omega _S$ is described in Kashiwara-Vergne [KV]. Here we restrict to $n = p$ and
$m=1$ and follow [LV] p.209f. \\

{\bf 4.6.1} From 3.3.2 we repeat the formulae for the derived \re of $\omega _S$ on $\Ss(\R^n)$ restricted to $\g'_c = \sl(2,\C) = <Z,X_{\pm}>$
$$
\hat Z = \pi(x,x) - (1/4\pi)\Delta ,\,\,\hat X_{\pm} = (1/2)( E + n/2 \mp (\pi(x,x) + (1/4\pi)\Delta)),\,\, E := \sum_{j=1}^n x_j\partial _j,
$$
and from 3.3.3 for the restriction to $\g = \{ Y \in \M(p,\R); Y = - \ta Y \}$
$$
\hat Y_{\alpha \beta } =  x_\alpha \partial _\beta + x_\beta \partial _\alpha .
$$
 As $\tilde G$ and $G'$ commute in $\tilde {\hat G}$, the subspace $L^2(\delta _m)$ of functions in $L^2(V) \simeq L^2(\R^n)$ of type 
 $\delta _m$  with respect to $\Or(S)$ is stable by $\tilde G \times \Or(S)$ and one has the multiplicity free decomposition
$$
L^2(V) = \oplus _{m \in \N_0} L^2(\delta _m).
$$ 
From 4.5.1 we take a harmonic polynomial of degree $m$ and put
$$
\varphi _P(x) := P(x) \varphi _0(x), \,\, \varphi _0(x) := e^{- \pi \ta xx} = e^{-\pi S(x)}.
$$
After a short calculation we get\\

{\bf 4.6.2} \Prop: 
$$
\hat Z \varphi _P = (m + p/2)\varphi _P,\,\,\hat X_- \varphi _P = 0.
$$
Hence $\varphi _P \in L^2(\delta _m)$ is a lowest weight vector for the discrete series \re $\pi^+_{m+p/2} =: \alpha $ 
of $G' =\SL(2,\R)$ resp.~its covering group $\GG$ and we have the decomposition of $\tilde \omega _S$ in unitary irreducible \res of 
$\Or(p) \times \tilde {\SL}(2,\R)$ as
$$
\tilde \omega _S = \oplus _m (\bar \pi_{(p/2) + m)} \otimes \delta _m).
$$ 

Following the same program as we did in 4.3, now we look at
$$
(\omega _S(g_\tau )\varphi _P)(x) = v^{\alpha /2}P(x)e^{\pi i \tau S(x)}
$$
and
$$ 
 \varphi _{P,\tau }(x) := I_\alpha (\omega _S(g_\tau )\varphi _P(x) = P(x)e^{\pi i \tau S(x)}.
$$
We have (at least for integral $\alpha$) the relation
$$
\omega _S(g)\varphi _{P,\tau } = (c\tau + d)^{-\alpha }\varphi _{P,g(\tau )};\,\, g = \begin{pmatrix} a&b\\c&d \end{pmatrix}
$$
and applying the theta distribution we get for the positive quadratic form $S$ the theta function with harmonic coefficient $P$
$$
\theta _P(\tau ) := \sum_{x \in \Z^p}  P(x)e^{\pi i \tau S(x)}
$$
for which one may expect some behavior of a modular form of weight $\alpha .$ \\

{\bf 4.6.3} As we did in 4.4, here we also can construct a Jacobi theta function living on the Jacobi group adapted to this case
$$
G^J_n := \SL(2,\R) \ltimes \Heis(\R^n).
$$
Analogously extending the previous notation, we get
$$
\omega_S ((p,q,\kappa )g_\tau r(\vartheta ))\varphi _P(x) = e^{\pi i(\kappa + 2(\ta x + \ta p)q + \tau S(x + p))} P(x+p)v^{\alpha /2}e^{i\alpha }
$$
and hence
$$
\omega_S ((p,q,\kappa )g_\tau r(\vartheta ))\varphi _P(0) = e^{\pi i(\kappa + 2\ta pq + \tau S(p))} P(p)v^{\alpha /2}e^{i\alpha } =: \Phi _P(g^J_n)
$$
a lowest weight vector living on $G^J_n$ which as in 4.4.5 may be used to introduce
$$
\begin{array}{rcl}
\Theta _P(g^J_n) & := & \sum_{\ell \in \Z^n} \Phi _P((\ell,0,0)g^J_n)\\*[0.3cm]
& = & e^{\pi i(\kappa + \ta p(\tau Sp + q))}v^{\alpha /2}e^{i\alpha \vartheta }\sum_{\ell \in \Z^p}P(p+\ell)e^{\pi i(\tau S(\ell) + 2\ta \ell(\tau Sp+q)} 
\end{array}
$$
For $P = 1$ with $\hat \tau _S := \tau S, \,\hat z_S := \tau Sp + q$ we see that this $\Theta _1$ is up to an automorphic 
lifting factor the special value (already coming up in 1.3.7) of the $\dim n-$Jacobi theta function
$$
\theta (\hat \tau ,\hat z) = \sum_{\ell \in \Z^n} e^{\pi i (\ta \ell \hat \tau \ell + 2\ta \ell\hat z)},
$$
which is discussed in [MuIII] Section 8.
It is interesting to analyze what happens if $P$ is a polynomial of higher degree. For this here we refer to [MuIII] Section 9 
treating Jacobi theta functions in spherical harmonics.\\

{\bf 4.6.4} We go back to 4.6.2 and for later use we reproduce from [LV] p.222 the introduction of an operator $\bar {\FF}_m$ which 
intertwines the \re $\tilde \omega _S$ of $\Or(S) \times \GG$ on $L^2(\delta _m)$ with the \re $\bar \pi_\alpha  \otimes \delta _m$ on the space 
$\bar {\OO}(\Ha) \otimes \Hi^m \simeq \bar {\OO}(\Ha,\Hi^m).$ This $\bar {\FF}_m$ is defined by
$$
<\bar {\FF}_m\psi  (\tau ), Q> := \int_V e^{- \pi i \bar \tau S(\xi )}\psi (\xi )Q(\xi )d\xi \fa \psi \in L^2(\delta _m),\,\, Q \in \Hi^m.
$$
$\bar {\FF}_m$ is injective and has the property
$$
\bar {\FF}_m \varphi _{P,\tau }(w) = c(\bar w - \tau )^{-\alpha } \otimes P.
$$

\subsection[Weil Representation Associated to an Indefinite Form]
{Decomposition of the Weil Representation Associated to an Indefinite Quadratic Form}

We take $G = \Or(V) \simeq \Or(p,q), G' = \SL(2,\R)$ and the covering $\GG.$ 
Adapting the notation to the one used here, we reproduce from [LV] 2.5.26\\

{\bf 4.7.1} \Th: The discrete spectrum of the \re $\tilde \omega _S$ of $G \times \GG$ is given as follows\\
A) For $p>1, q>1$
$$
(\tilde \omega _S)_d = \oplus_{\alpha >1}(\hat \delta _\alpha \otimes \bar \pi _\alpha ) \oplus \oplus _{\beta >1}(\hat \delta _\beta \otimes \pi _\beta )
$$
where $\alpha ,\beta \in \Z$ if $(p-q)/2 \in \Z$ and $\alpha ,\beta \in (1/2)\Z$ if $(p-q)/2 \in (1/2)\Z.$ 
The \re $\hat \delta _\alpha $(resp.~$\hat \delta _\beta $) is a irreducible \re of $\Or(p,q).$ Its restriction to $\Or(p) \times \Or(q)$ is
$$
\hat \delta _\alpha = \oplus _{k,m} \delta _k \otimes \delta _m, \,\, k-m+(p-q)/2 = \alpha + 2j, \,\, j\geq 0
$$
resp.
$$
\hat \delta _\beta = \oplus _{k,m} \delta _k \otimes \delta _m, \,\, m-k+(q-p)/2 = \beta + 2j, \,\, j\geq 0.
$$
B) For $p>1, q=1$
$$
(\tilde \omega _S)_d = \oplus_{\alpha >1}(\hat \delta _\alpha \otimes \bar \pi _\alpha )
$$
with
$$
\hat \delta _\alpha = \oplus _{k,m=0,1} \delta _k \otimes \delta _m, \,\, k-m+(p-1)/2 = \alpha + 2j, \,\, j\geq 0.
$$
C) For $p=q=1$
$$
(\tilde \omega _S)_d = 0.
$$
We will discuss this theorem and show a way to get lowest weight vectors for these \res. 
The first thing to remark is that there is a kind of precursor to this theorem  going back to Gutkin and Repka 
treating the decomposition of the tensor product of two discrete series \res of $\SL(2,\R).$\\

{\bf 4.7.2} \Rem: The product $\bar \pi_\alpha  \otimes \pi_\beta , \alpha \geq \beta $ contains discretely the sum
$$
\oplus _{\alpha - \beta - 2j>1, j\in\Z} \,\,\bar \pi_{\alpha - \beta - 2j}.
$$
{\bf 4.7.3} Next, we write according to the orthogonal decomposition $V = V_1 \oplus V_2,\,V_1 \simeq \R^p, V_2 \simeq \R^q$
$$
L^2(V) = L^2(V_1) \otimes L^2(V_2).
$$
Then from 4.7.1 we know that $\tilde \omega _S$ as a \re of \,$\tilde {\SL}(2,\R) \times \Or(p) \times \Or(q)$ is isomorphic to 
$$
\oplus _{k,m}\,\, \bar \pi_{(p/2)+k} \otimes \pi_{(q/2)+m} \otimes \delta _k \otimes \delta _m.
$$
 We put
$$
d := k + p/2 - (m + q/2)
$$
and will see that, when $d > 1,$ the \re $\bar \pi_{k + p/2} \otimes \pi_{m + q/2}$ contains $\bar \pi_d$ by analyzing 
the shape of a lowest weight vector $v_0$ of $\omega _S:$\\

{\bf 4.7.4} From 3.3.2 we recall the formulae for $\s(2,\R)_c$
$$
\begin{array}{rcl}
\hat Z & = & \pi(x,x) - (1/(4\pi))\Delta \\
\hat X_{\pm} & = & (1/2)(E + n/2 \mp (\pi(x,x) + (1/(4\pi))\Delta ))
\end{array}
$$
where
$$
(x,x) = S(x) := \sum_{\alpha =  1}^p x^2_\alpha - \sum_{\mu = p+1}^{p+q} x_\mu ^2, \,\,
\Delta = \sum_{\alpha =  1}^p \partial _{x_\alpha}^2 - \sum_{\mu = p+1}^{p+q} \partial _{x_\mu} ^2,\,\,E = \sum_{j=1}^{p+q} x_j\partial _{x_j}.
$$
and from 3.3.3 for $\ot(p,q)$ with  $1\leq \alpha ,\beta \leq p,\,\,p+1\leq \mu ,\nu \leq p+q$
$$
\hat Y_{\alpha \beta } = x_\alpha\partial _{x_\beta } - x_\beta \partial _{x_\alpha },\,\, 
\hat Y_{\mu \nu } = x_\mu \partial _{x_\nu } - x_\nu \partial _{x_\mu }, \,\,
\hat Y_{\alpha \mu } = x_\alpha\partial _{x_\mu } + x_\mu \partial _{x_\alpha }.
$$
Our $\tilde \omega _S$ on $\Hi = L^2(\R^n)$ is simultaneously a \re of $G = \Or(p,q)$ and $\tilde G' = \Mp(2,\R).$ 
Hence it is natural to look at subspaces of functions which are invariant under both groups and 
we take
$$
\Hi_+ := \{ \varphi \in \Hi ; \varphi \vert_{\X^-} = 0 \}
$$
We recall the remarks 3.3.5 and 3.4.1 and look at functions
$$
\varphi = \psi \varphi _1, \,\, \varphi _1(x) = e^{-\pi S(x)},\,\,\psi \in \Cl^\infty (\R^n)
$$
for $x$ with $S(x) \geq 0$ and $\varphi = 0$ for $S(x) < 0.$ By a small calculation we get\\

{\bf 4.7.5} \Rem: One has
$$
\hat Z \varphi = (-(1/4\pi))\Delta \psi + E\psi + (n/2)\psi) \varphi _1,\,\, \hat X_- \varphi = (1/(8\pi)) \Delta \psi \varphi _1.
$$
We see that $\varphi = \psi e^{-\pi S(x)}$ is a lowest weight vector of weight $\lambda $ if 
one has
$$
\Delta \psi = 0,\,\, E\psi + (n/2)\psi = \lambda \psi,\,\, \varphi = \psi e^{-\pi S(x)} \in \Hi.
$$
Hence, we have to look for $\psi$ fulfilling these conditions.
If $\psi$ is a homogeneous polynomial of degree $m$ 
which is annihilated by the Laplacian $\Delta $ belonging to the indefinite form $S$, this could lead to the weight $m+n/2$. 
But that's not the true story. 
A refined discussion of the solution of these equations and the \res showing up in the decomposition is done 
by Rallis and Schiffmann in [RS3] with a summary of their results given in [RS1] or [RS2]. As already done above, here we follow the 
version given by Vergne in [LV] p.225f. \\

{\bf 4.7.6} \Rem: By some easy calculation one has with $\Delta , E,$ and $S$ as in 4.7.5
$$
\Delta (S^\alpha f) = S^\alpha \Delta f + 4\alpha S^{\alpha - 1}(E + (n/2) + \alpha - 1)f, \,\, f \in \Cl^\infty (\R^n)
$$
and with\\*[0.2cm]
\phantom{}\hspace{1cm} $P_1$ a harmonic polynomial of degree $k$ in $x_1,\dots,x_p,$ \\
\phantom{}\hspace{1cm} $S_1 = \sum_{\alpha = 1}^p x_\alpha ^2,$\\
\phantom{}\hspace{1cm} $P_2$ a harmonic polynomial of degree $m$ in $x_{p+1},\dots,x_{p+q},$ \\
\phantom{}\hspace{1cm} $S_2 = \sum_{\mu = 1+p}^{p+q} x_\mu ^2,$\\
and\\
\phantom{}\hspace{1cm} $\psi  := P_1P_2 S_1^\alpha S_2^\beta  S^\gamma$\\ 

for $\,S_1(x) \not=0$ and $S_2(x) \not=0\,$ one has
$$
\Delta \psi = 0
$$
if\\
\phantom{}\hspace{1cm} 1) $(p+q)/2 + \gamma  - 1 + 2\alpha + 2\beta + k + m = 0$ and\\
\phantom{}\hspace{1cm} 2) $\alpha (p/2 + \alpha  - 1 + k) = 0$ and \\
\phantom{}\hspace{1cm} 3) $\beta (q/2 + \beta  - 1 + m) = 0.$\\

If we put $\beta = 0$ we already get the first statement in the following theorem.\\

{\bf 4.7.7} \Th: For 
$$
\begin{array}{rcll}
\psi _{P_1,P_2} & := & P_1P_2 S_1^{k + (p-2)/2}S^{(p-q)/2 + k - m - 1} & \fa x \with S(x) > 0\\
& = & 0  & \fa x \with S(x) < 0 
\end{array}
$$
and the {\it Rallis-Schiffmann function}
$$
 \varphi _{P_1,P_2} := \psi _{P_1,P_2}e^{- \pi S(x)}
$$
one has with $d := k - m + (p - q)/2$\\

\phantom{}\hspace{1cm} 1) $\Delta \psi _{P_1,P_2} = 0,\,\, (E + n/2) \psi _{P_1,P_2} = d \cdot \psi _{P_1,P_2}$\\
\phantom{}\hspace{1cm} 2) $\varphi _{P_1,P_2} \in L^2(\R^n)$ if $d > 1,$\\
\phantom{}\hspace{1cm} 3) $\varphi _{P_1,P_2} \in L^1(\R^n)$ if $k-m > q,$\\
\phantom{}\hspace{1cm} 4) If $p+q > 2$ and $k-m > q,$ then $\varphi _{P_1,P_2}$ is in $L^1(\R^n) \cap L^2(\R^n)$ and is continuous.\\

As already said, item 1) is clear from our previous remarks. For the other statements we refer to [LV] p.228/9.\\

{\bf 4.7.8} \Cor: Hence the Rallis-Schiffmann function is a lowest weight vector of weight $d$ for the \re of $\tilde G'$ and of type $k$ and $m$ with 
respect to the action of $\Or(p) \times \Or(q).$\\

 Vergne observes that moreover for the indefinite case one has an analogy to the definite formulae in 4.6.2. 
Using the operator  $\FF_m$ introduced in 4.6.4 one sees that the operator $\bar {\FF}_k \otimes \FF_m$ intertwines the \re 
$\tilde \omega _S$ of 
$\tilde G' \times \Or(p) \times \Or(q)$ restricted to $L^2(\delta _k \otimes \delta _m)$ with the \re
$$
\bar \pi_{(p/2)+k} \otimes \pi_{(q/2)+m} \otimes \delta _k \otimes \delta _m.
$$
The \re $\bar \pi_{(p/2)+k} \otimes \pi_{(q/2)+m}$ acts on the space of functions $F(\bar \tau _1,\tau _2)$ antiholomorphic in $\tau _1$ and 
holomorphic in $\tau _2$ by
$$
((\bar \pi_{(p/2)+k} \otimes \pi_{(q/2)+m})(g^{-1})F)(\bar \tau _1,\tau _2) = 
 (c\bar \tau _1+d)^{-((p/2)+k)}(c\tau +d)^{-((q/2)+m)}F(g(\bar \tau _1),g(\tau _2)).
$$
A function antiholomorphic in one and holomorphic in the other variable is is entirely determined by its restriction to the diagonal $(\bar \tau ,\tau ).$ 
Hence it is natural to consider the \re $\bar \pi_{\alpha ,\beta }$ of $\tilde{ \SL}(2,\R)$ acting on all functions on $\Ha$ by
$$
(\bar \pi_{\alpha ,\beta }(g^{-1})f)(w) :=  (c\bar \tau _1+d)^{-\alpha )}(c\tau +d)^{-\beta }f(\tau (g)),\,\, w \in \Ha.
$$
The operator 
$$
\varphi \longmapsto ((\bar {\FF}_k \otimes \FF_m)\varphi )(\bar \tau ,\tau )
$$
intertwines the \re $\tilde \omega _S\vert_{L^2(\delta _k\otimes \delta _m)}$ with $\bar \pi_{(p/2)+k,(q/2)+m} \otimes \delta _k \otimes \delta _m.$ 
We still denote this operator as $(\bar {\FF}_k \otimes \FF_m)$ and introduce another operator by
$$
f(w) \longmapsto (Mf)(w) := ({\rm Im} \,w)^{-((q/2)+m)}f(w),
$$
which intertwines the \re $\bar \pi_{d,0},\,\, d = ((p/2) + k - ((q/2) + m)$ 
acting on the functions on $\Ha$ with the \re $\bar \pi_{(p/2)+k,q/2)+m}.$ 
The \re $\bar \pi_d$ is naturally contained in $\bar \pi_{d,0},$ thus in $\bar \pi_{(p/2)+k,(q/2)+m}.$ 
From 4.2.2 we deduce that the function on $\Ha$ given by 
$$
\psi '_\tau (w) = ({\rm Im} \,w)^{-((q/2)+m)}(\bar w - \tau )^{-d}
$$
verifies the relation (excuse the double meaning of the letter $d$)
$$
\bar \pi_{(p/2)+k,(q/2)+m}(g)\psi '_\tau  = (c\tau  + d)^{-d}\psi '_{g(\tau )}.
$$ 
As to be seen by some calculation as in [LV] p.234 the variant of the Rallis-Schiffmann function from 4.7.7
$$
\varphi _{P_1P_2,\tau} (x) := P_1P_2S_1^{-(k+(p-2)/2)}S^{d-1}e^{\pi i\tau S(x)},\,\, \tau \in \Ha
$$
 has the following properties.\\

{\bf 4.7.9} \Th: For $d > 1$ and harmonic polynomials $P_1,P_2$ of degree $k$ resp.~$m$ in $p$ resp.~$q$ variables one has
$$
 \varphi _{P_1P_2,\tau }(x) \in L^2(\delta _k \otimes \delta _m)
$$
and 
$$
((\bar {\FF}_k \otimes \FF_m)  \varphi _{P_1P_2\tau }(x) = \psi' _\tau \otimes P_1 \otimes P_2.
$$
Hence $ \varphi _{P_1P_2,\tau }(x)$ fulfills the fundamental formula
$$
\omega _S(g)  \varphi _{P_1P_2,\tau} (x) = j(g,\tau )^{-d} \varphi _{P_1P_2g(\tau)} (x)
$$
which makes it a candidate for the production of a modular form of weight $d$ via a theta distribution.\\

{\bf 4.7.10} \Rem: We remind that there is a direct way from the lowest weight vector $\varphi _{P_1P_2}$ to $\varphi _{P_1P_2\tau }(x).$ 
From the formulae of the Weil \re with our matrix $g_\tau $ transforming $i$ to $\tau \in \Ha$ as already used in 4.3.2 we have in this  case
$$
\omega (g_\tau )\varphi _{P_1P_2}(x) = v^d \varphi _{P_1P_2g(\tau)} (x).
$$

{\bf 4.7.11} Up to now we treated the Rallis-Schiffmann functions only with respect to their behaviour 
concerning the group $\tilde {\SL}(2,\R).$ But one can proceed similarly concerning the orthogonal group. 
Here we take as an example the case $p = 2, q = 1.$ In 3.3.4 we determined operators for the complexified algebra
$$
\hat H_0 = - 2i(x_1\partial _2 - x_2\partial _1),\,\,\hat Y_{\pm} = - (x_1 \pm ix_2)\partial _3 - x_3(\partial _1 \pm i\partial _2).
$$
For 
$$
\varphi = \psi \varphi _1,\,\, \varphi _1(x) = e^{-S(x)}
$$
one gets
$$
\hat H_0 \varphi  = -2i(x_1\psi _{x_2} - x_2\psi _{x_1})\varphi _1,\,\,
\hat Y_{\pm}\varphi = - (x_1\psi _{x_3} + x_3\psi _{x_1} \mp i(x_2\psi _{x_3} + x_3\psi _{x_2})\varphi _1.
$$
For the Rallis-Schiffmann function from Theorem 4.7.7 in this case one has only  two choices for the polynomial $P_2,$ namely 
$P_2(x) = 1$ or $= x.$ We take $P_2 = 1$ and get
$$
 \varphi _{P_1,P_2} := \psi _{P_1,P_2}e^{- \pi S(x)}
$$
with
$$\begin{array}{crll}
\psi _{P_1,P_2} & := & P_1 S_1^{k }S^{k - 1/2} & \fa x \with S(x) > 0\\
& = & 0  & \fa x \with S(x) < 0 
\end{array}
$$
As a homogeneous harmonic polynomial $P_1$ in two variables one can take $P_1(x_1,x_2) := (x_1 \pm ix_2)^k.$ Then we have
$$
\varphi _{\pm}^k := \psi^k _{\pm}\varphi _1,\,\, \psi ^k _{\pm} := (x_1 \mp ix_2)^{-k}S(x)^{k-1/2}
$$
and by the formula above get
$$
\hat H_0 \varphi ^k _{\pm} = \pm 2k \varphi ^k _{\pm};\,\, \hat Y_- \varphi ^k _{+} = 0, \,\, \hat Y_+ \varphi ^k _{-} = 0.
$$
Hence we get a refinement of corollary 4.7.8.\\

{\bf 4.7.12} \Prop: $\varphi ^k _{+}$ is simultaneously a lowest weight vector of weight $2k$ for a \re $\hat \delta _\alpha , \alpha = k+1/2$ of 
$G = \Or(2,1)$ and of weight $k + 1/2$ for the \re $\bar \pi_\alpha $ of $\tilde G' = \Mp(2,\R).$\\

In this example we finally have the corner stone for the construction of a theta function living simultaneously 
on the orthogonal and the metaplectic group and, hence, apt to produce a correspondence between automorphic forms 
belonging to these groups.\\

{\bf 4.7.13} We pursue this a bit reproducing the construction of the theta function leading to the Shimura correspondence from [LV] p.268f.\\

As orthogonal space we consider the vector space 
$$
E := \Sym_2(\R) = \{ y =  \begin{pmatrix} y_1&y_3\\y_3&y_2 \end{pmatrix}; \,y_1,y_2,y_3 \in \R \}
$$
with the quadratic form $S'(y) = - 2\det y$ and the associated bilinear form
$$
S'(y,y') = 2y_3y'_3 - y_1y'_2 - y_2y'_1.
$$
The transformation
$$
y_1 = x_3 + x_1, \,\, y_2 = x_3 - x_1,\,\,y_3 = x_2
$$
leads to our usual signature $(2,1)$ situation
$$
S'(y) = 2(x_1^2 + x_2^2 - x_3^2).
$$
The group $\SL(2,\R)$ acts on $E$ by $g\cdot y := gy\ta g.$ This action leaves $S(y)$ stable and - as well known - 
leads to a surjective map $\phi : \SL(2,\R) \longrightarrow \Or(2,1)^0$ with $\ker \phi = \{ \pm 1_2 \}.$ 
This map is given by
$$
 g = \begin{pmatrix} a&b\\c&d \end{pmatrix} \longmapsto  \phi (g) := \begin{pmatrix} a^2&b^2&2ab\\c^2&d^2&2cd\\ac&bd&ad+bc \end{pmatrix}.
$$
We use $\phi $ to transfer the restriction $\tilde \omega _S$ of the Weil \re from $\Or(2,1) \times \tilde G'$ to a \re of 
$\SL(2,\R) \times \tilde G'$ where the \re of the first group, abbreviated by $\tilde G,$ is given by 
$\varphi (y) \longmapsto \varphi (\phi (g^{-1})y).$ Our Rallis-Schiffmann function written in the $y-$coordinates looks like
$$
\varphi _+(y) = ((1/2)(y_1 - y_2) - iy_3)^{-k}S'(y)e^{-\pi S'(y)}.
$$
For $z = x +iy$ and $g = g_z =  \begin{pmatrix} y^{1/2}&xy^{-1/2}\\&y^{-1/2} \end{pmatrix}$ one has
$$
\phi (g_z^{-1}) =  \begin{pmatrix} y^{-1}&x^2y^{-1}&-2xy^{-1}\\0&y&0\\0&-x&1 \end{pmatrix}
$$
and $\phi (g_z^{-1})$ transforms 
$$
((1/2)(y_1 - y_2) - iy_3 \longmapsto (1/(2y))(y_1 + z^2y_2 - 2zy_3) = - S'(y,Q(z)); \,\, Q(z) = \begin{pmatrix} z^2&z\\z&1 \end{pmatrix}.
$$ 
{\bf 4.7.14} Hence, using 4.7.10, we get by application of the Weil \re to our Rallis-Schiffmann function expressed in the $y-$coordinates
$$
(\omega (g_z \cdot g'_\tau )\varphi _{P_1P_2})(y) = (-2y)^kv^{k+(1/2)/2}S'(y,Q(z)) S'(y)^{k-1/2}e^{\pi i \tau S'(y)}.
$$
for $S'(y) > 0$ and by zero for $S'(y) \leq  0.$ We put
$$
\varphi ^k(z,\tau )(y) := S'(y,Q(z)) S'(y)^{k-1/2}e^{\pi i \tau S'(y)}
$$
and have a function with the fundamental relation
$$
\tilde \omega (g,g')\varphi ^k(z,\tau ) = j(g,z)^{-2k}j(g',\tau )^{-(k+1/2)}\varphi ^k(g(z),g'(\tau )).
$$
By a suitable skillful averaging procedure (as in [LV] p.272f) one comes to a theta function in both variables $z,\tau :$ 
We consider the groups
$$
\Gamma _0(N) := \{ \begin{pmatrix} a&b\\c&d \end{pmatrix} ; \,\, c \equiv 0 \mod N \}.
$$
and
$$
\Gamma _0(0,2N) := \{ \begin{pmatrix} a&b\\c&d \end{pmatrix} ; \,\, b \equiv 0 \mod 2N \}.
$$
$\tilde \Gamma _0(N) $ denotes the inverse image of $\Gamma _0(N) $ in $\tilde {\SL}(2,\R).$

Let $\psi $ be character mod $4N$ and with a slight abuse of notation also the character
$$
\gamma = \begin{pmatrix} a&b\\c&d \end{pmatrix} \longmapsto \psi (d)
$$
of $\Gamma _0(4N).$ $\lambda $ denotes the theta character, i.e.~for $\tilde \gamma \in \tilde \Gamma _0(4)$
$$
\lambda (\tilde \gamma ) := \varepsilon _d(\frac{c}{d}), 
$$
and $u$ a function on $\Z/N\Z$ satisfying $u(aj) = \psi (a)u(j).$ The one has the central statement.\\

{\bf 4.7.15} \Th (Theorem 2.7.17 in [LV]): For $k>1$ the function $\Omega _u$ given by
$$
\Omega _u (z,\tau ) := \sum _{y \in \Z^3; S'(y) > 0} u(y_1)\varphi ^k(z,\tau )(y)
$$
is a holomorphic function of $(z,\tau )$, which is\\

\phantom{} \hspace{.5cm}- modular in $\tau $ with respect to $\tilde {\Gamma} _0(4N)$, with character $\lambda \psi ,$ of weight $k + (1/2)$\\
and\\
\phantom{} \hspace{.5cm}- modular in $z$ with respect to $\Gamma _0(0,2N),$ with character $\psi^{-2},$ of weight $2k.$\\

For a function $u$ on $\Z/N\Z$ the Fourier transform $\hat u$ is defined by 
$$
\hat u(m) := \sum_{h \in \Z/N\Z} u(h) e^{- 2\pi i mh/N}.
$$
If one chooses $u = u_0$ with $\hat u_0 = \bar \psi $ the Petersson inner product of a certain cusp form of weight $k+(1/2)$ with 
$\Omega _{u_0}$ produces an automorphic form of weight $2k$ and thus establishes a version of the {\it Shimura correspondence.} 
This is only an example of much more material which has been obtained in a similar fashion. We refer to the other 
chapters of [LV] and, for instance to [RS2].\\

{\bf 4.7.16} There is the challenge to try to extend parts of this to pairs consisting of an orthogonal group and a Jacobi group or 
even an euclidean group and a Jacobi group. But this is no longer prehistory and has to appear again later.\\

The following list contains some items (but not all) which don't belong to Part I and will be needed in the other Parts of this text.\\


\printindex
Rolf Berndt\\
Mathematisches Seminar der Universit\"at Hamburg\\
Bundesstr. 55\\
D-20146 Hamburg\\
Germany\\
berndt@math.uni-hamburg.de

\end{document}